\documentclass{amsart}
\usepackage{amsmath}%
\usepackage{amsfonts}%
\usepackage{amssymb}%
\usepackage{graphicx}
\usepackage{mathtools, tikz,float}

\newcommand{\p}{\mathcal{P}}
\renewcommand{\P}{\mathbb{P}}
\newcommand{\F}{\mathcal{F}}
\newcommand{\G}{\mathcal{G}}
\newcommand{\C}{\mathcal{C}}
\renewcommand{\L}{\mathcal{L}}
\newcommand{\M}{\mathcal{M}}

\newcommand{\D}{\mathcal{D}}

\DeclareMathOperator{\Mod}{Mod}

\newcommand{\R}{\mathbb{R}}

\newcommand{\h}{\mathbb{H}^2}
\renewcommand{\H}{\mathcal{H}}
\renewcommand{\S}{\mathcal{S}}
\newcommand{\T}{\mathcal{T}}

\newtheorem{theorem}{Theorem}[section]
\theoremstyle{plain}
\newtheorem{lem}[theorem]{Lemma}
\newtheorem{claim}[theorem]{Claim}
\newtheorem{cor}[theorem]{Corollary}
\newtheorem{rem}[theorem]{Remark}
\newtheorem{prop}[theorem]{Proposition}

\newtheorem{ques}{Question}

\numberwithin{equation}{subsection}
\theoremstyle{definition}
\newtheorem{defi}[theorem]{Definition}

\title{Bounds on the number of non-simple closed geodesics on a surface}
\author{Jenya Sapir}

\begin{document}
\begin{abstract}
We give bounds on the number of non-simple closed curves on a
negatively curved surface, given upper bounds on both length and self-intersection
number. In particular, it was previously known that the number of all closed curves of length at most $L$ grows exponentially in $L$. We get exponentially tighter bounds given weak conditions on self-intersection number.
\end{abstract}

\maketitle
\section{Introduction}
Let $\S$ be a genus $g$ surface with $n$ boundary components, and let $X$ be a negatively curved metric on $\S$.  Closed geodesics on surfaces have been studied extensively over the years. In this paper, we give upper and lower bounds on the number of closed geodesics on $\S$ given upper bounds on length and self-intersection number. The lower bound follows from a lower bound on the number of closed geodesics in a pair of pants, and is proven in our earlier paper, \cite{Sapir}. The upper bound comes from looking at closed geodesics on a closed surface with a flat metric.

\subsection{Statement of results}
Let $\G^c$ be the set of all closed geodesics on a surface $\S$ with negatively curved metric $X$. Then let 
\[
 \G^c(L,K) = \{\gamma \in \G^c \ | \ l_X(\gamma) \leq L, i(\gamma, \gamma) \leq K\}
\]
where $l_X(\gamma)$ is the length of $\gamma$ in $X$, and $i(\gamma, \gamma)$ is its geometric self-intersection number. We will write $\#\G^c(L,K)$ for the number of curves in $\G^c(L,K)$. We wish to get bounds on $\#\G^c(L,K)$ in terms of both $L$ and $K$. As a first step, we prove the following bounds when $K$ is fixed:
\begin{theorem}
\label{thm:AppendixIntro}
 Fix $K \geq 0$. Let $\S$ be a closed surface of genus $g$, with hyperbolic metric $X$. Then
 \[
 \#\G^c(L,K) \asymp L^{6g-6}
 \]
where the constants depend only on the constant $K$ and the metric $X$.
\end{theorem}
Note that we write $A(L) \asymp B(L)$ if there are constants $a$ and $b$ so that $\frac 1 a B(L) - b \leq A(L) \leq a B(L) + b$. This means that the number of curves with a finite bound on self-intersection number grows at the same rate as the number of simple closed curves. 

However, the dependence on the constant $K$ in this theorem is not explicit. We also give upper and lower bounds with explicit dependence on $K$. The upper bound is as follows:
\begin{theorem}
\label{thm:MainUpperBound}
For any negatively curved metric $X$ on $\S$, and for any $L \geq 0, K \geq 1$, we have
 \[
 \# \G^c(L,K) \leq \min 
 \left \{  c_X e^{\delta L} ,  (c_X L)^{c_\S \sqrt K} \right \}
\]
where $c_X$ depends on $X$, $c_\S$ depends only on $\S$, and $\delta$ is the topological entropy of the geodesic flow on $\S$ with respect to $X$.

\end{theorem}

Margulis \cite{MargulisPhD} gave the asymptotic growth of the size of $\G^c(L)$, the set of all closed geodesics of length at most $L$ (see below). This gives us an upper bound on $\#\G^c(L,K)$ for any $K$. Theorem \ref{thm:MainUpperBound} gives an exponentially better upper bound whenever $K(L) = o (\frac{L^2}{\ln^2 L})$:
\begin{cor}
 If $K = K(L)$ is a function of $L$ such that $K = o (\frac{L^2}{\ln^2 L})$, then for any $0 < c < 1$,
 \[
  \frac{\# \G^c(L,K)}{\# \G^c(L)} < e^{-cL}
 \]
 for all $L$ large enough, depending on $c$ and $X$.
\end{cor}

The following lower bound on $\#\G^c(L,K)$ is proven in our earlier paper \cite{Sapir}:
\begin{theorem}
\label{thm:LowerBoundSurface}
 Let $X$ be a hyperbolic metric on $\S$. Then whenever $K > 12$ and $L > 6 s_X \sqrt K$ we have 
 \[
 \# \G^c(L,K) \geq c_X \left (\frac{L}{6 \sqrt K} \right )^{6g-6+2n}2^{\sqrt{\frac{ K}{ 12}}}
 \]
 where $s_X$ and $c_X$ are constants that depend only on the metric $X$.
\end{theorem}

As $L$ goes to infinity, this theorem suggests a way to interpolate between the case when $K$ is a constant and the case when $K$ grows like $L^2$. If $K$ is a constant, and $L$ is  large enough, this theorem says $\#\G(L,K) \geq c'_X L^{6g-6+2n}$, for $c'_X$ a new constant. This is consistent with Theorem \ref{thm:AppendixIntro}. For $K = O(L^2)$, however, we have that $\frac{L}{3 \sqrt K} = O(1)$, and Theorem \ref{thm:LowerBoundSurface} gives an exponential lower bound on $\#\G^c(L,K)$ in $L$ that is consistent with exponential growth for the set of all closed geodesics. This theorem demonstrates the transition from polynomial to exponential growth of the number of geodesics on $\S$ in terms of their length and self-intersection number.

\subsection{Previous results on an arbitrary surface}
The problem of counting closed geodesics in many contexts has been studied extensively. There is an excellent survey of the history of this problem by Richard Sharp that was the  published in conjunction with Margulis's thesis in \cite{Margulis04}. 

In brief, let $\G^c$ be the set of closed geodesics on $\S$ and let 
\[
 \G^c(L) = \{\gamma \in \G^c \ | \ l(\gamma) \leq L\}
\]
where $l(\gamma)$ is the length of $\gamma$. Then Margulis showed that for a finite area, negatively curved surface,
\begin{equation}
\label{GL}
 \#\G^c(L) \sim \frac{e^{\delta L}}{\delta L}
\end{equation}

The number of closed geodesics with upper bounds on length have since been thoroughly studied. As a next step, we can count closed geodesics with respect to both length and self-intersection number. 

Bounds on the number of simple closed curves were first given in \cite{Rees81}. Mirzakhani then showed that for a hyperbolic surface $\S$ of genus $g$ with $n$ punctures,
\[
 \#\G^c(L,0) \sim c(\S) L^{6g-6 + 2n}
\]
where $c(\S)$ is a constant depending only on the geometry of $\S$ \cite{Mirzakhani08}. This result was extended by Rivin  to geodesics with at most one self-intersection, to get that 
\[
  \# \G^c(L,1) \sim c'(\S) L^{6g-6 +2n}
\]
where $c'(\S)$ is another constant depending only on the geometry of $\S$ \cite{Rivin12}. 

For arbitrary functions $K = K(L)$, no asymptotic bounds were known. We can instead ask the following question as a first step to finding asymptotics.
\begin{ques}
 Given arbitrary $L$ and $K$, what are the best upper and lower bounds we can get on $\#\G^c(L,K)$?
\end{ques}

Trivial bounds come from the fact that $\#\G^c(L,0) \leq \#\G^c(L,K) \leq \#\G^c(L)$,  but these bounds do not have any dependence on $K$. The theorem in the appendix gives a first bound for fixed $K$, but we do not get an explicit dependence on self-intersection number. Our main theorem gives bounds that are explicit in both length $L$ and intersection number $K$.

%

This paper is part of the author's PhD thesis, which was completed under her advisor, Maryam Mirzakhani. The author would especially like to thank her for the many conversations that led to this work. The author would also like to thank Jayadev Athreya, Steve Kerckhoff and Chris Leininger for their help and support. 

\section{Fixed intersection number}
Consider the family of sets $\{\G^c(L,K)\}_{L \geq 0}$, where $K$ is fixed but $L$ goes to infinity. We show that the size of these sets grows like a polynomial in $L$ of degree $6g-6$, which is the same as for simple closed curves.
\begin{theorem}
\label{thm:Appendix}
 Fix $K \geq 0$. Let $\S$ be a closed surface of genus $g$, with hyperbolic metric $X$. Then
 \[
 \#\G^c(L,K) \asymp L^{6g-6}
 \]
where the constants depend only on the constant $K$ and the metric $X$.
\end{theorem}

Note that we write
\[
 A(L) \asymp B(L) \iff \frac 1c B(L) - d \leq A(L) \leq c B(L) + d
\]
 for some constants $c$ and $d$ independent of $L$. We will write $A(L) \lesssim B(L)$ or $A(L) \gtrsim B(L)$ if only the left-hand or the right-hand inequality holds, respectively.
 
The idea for the proof is as follows. Let $\Mod_g$ denote the mapping class group of our genus $g$ surface $\S$. Then $ \Mod_g$ acts on $\G^c$. For each $f \in \Mod_g$, we let $f \cdot \gamma$ be the geodesic representative of $f$ applied to $\gamma$. Note that $i(\gamma, \gamma) = i(f \cdot \gamma, f \cdot \gamma)$. In other words, all curves in the $\Mod_g$ orbit of $\gamma$ have the same self-intersection number. Let $\Mod_g \cdot \gamma$ denote the $\Mod_g$ orbit of $\gamma$, and let
\[
 N(\gamma, L) = \# \{ f \cdot \gamma \in \Mod_g \cdot \gamma \ |\ l_X(f \cdot \gamma) \leq L\}
\]
where $l_Y(\gamma)$ denotes the length of $\gamma$ with respect to some metric $Y$ on $\S$. This is the number of curves in $\Mod_g \cdot \gamma$ of length at most $L$.
Then we can write $\#\G^c(L,K)$ as a sum:
\[
 \# \G^c(L,K) = \sum N(\gamma, L)
\]

There are finitely many $\Mod_g$ orbits of curves with at most $K$ self-intersections. To see this, imagine cutting $\S$ along $\gamma$. An Euler characteristic argument implies that this gives us at most $K$ connected components, each with piecewise geodesic boundary. Each connected component must fall into one of finitely many homeomorphism types. The total number of geodesic boundary arcs among all the components will be at most $4K$ (since $\gamma$ has at most $2K$ simple pieces between self-intersections). There are finitely many ways to choose at most $K$ (non-distinct) homeomorphism types of pieces, assign each piece at most $2K$ boundary arcs, and glue these shapes back together into $\S$. Therefore the above sum is finite, and we can bound $\#\G^c(L,K)$ by bounding $N(\gamma, L)$ for each $\gamma$.

The constants in Theorem \ref{thm:Appendix} depend on the number $f(K)$ of $\Mod_g$ orbits of closed curves with at most $K$ self-intersections. The number $f(K)$ is only known for finitely many $K$. Therefore, this theorem does not, in general, give an explicit dependence of $\#\G^c(L,K)$ on self-intersection number. See Section \ref{sec:Dependence} for more details on the dependence of the bounds in Theorem \ref{thm:Appendix} on $K$.

\subsection{Filling curves}

It is easier to bound the size of $N(\gamma, L)$ when $\gamma$ is filling, so we do this first.

\begin{lem}
\label{lem:filling}
 Let $\gamma$ be a filling curve. Then 
 \[
  N(\gamma, L)  \asymp  L^{6g-6}
 \]
where the constants depend only on $X$ and the orbit $\Mod_g \cdot \gamma$.
\end{lem}
\begin{proof}
We wish to count the number of elements in the $\Mod_g$ orbit of $\gamma$ whose length is at most $L$. Note that
\[
 l_X(f \cdot \gamma) = l_{f^{-1}X}(\gamma)
\]
In other words, instead of looking at the $X$-length of curves in the $\Mod_g$ orbit of $\gamma$, we can look at the $f X$-length of $\gamma$ for various $f \in \Mod_g$. Thus, 
 \[
 N(\gamma, L) =  \# \{ f \cdot \gamma  \in \Mod_g \cdot \gamma \ |\ l_{f^{-1}X}(\gamma) \leq L\}
 \]

We wish to restate this as a problem of counting $\Mod_g$ orbit points of $X$ in Teichmuller space. Let $\Mod_g \cdot X$ denote the $\Mod_g$ orbit of $X$. As an intermediate step, let
\[
 N_\gamma(X,L) = \#\{f X \in \Mod_g \cdot X \ | \ l_{f X}(\gamma) \leq L\}
\]
be the number of points in the $\Mod_g$ orbit of $X$ for which $\gamma$ has length at most $L$. Note that the set of $fX$ with $l_{fX}(\gamma) \leq L$ is the same as the set of $f^{-1}X$ with $l_{f^{-1}X}(\gamma) \leq L$. So we drop the inverses in what follows. 

We will first relate $N_\gamma(X,L)$ to $N(\gamma,L)$. The problem with the upper bound is those $f \in \Mod_g$ so that $f X = X$ but $f \cdot \gamma \neq \gamma$. By \cite{Hurwitz}, the automorphism group of $X$ has size at most $84(g-1)$. Thus,
\[
 N(\gamma, L) \leq 84(g-1)N_\gamma(X,L)
\]

Similarly, the stabilizer of $\gamma$ becomes an issue for the lower bound. Let 
\[
 \mathcal H = \{ f \in \Mod_g \ |\ f \cdot \gamma = \gamma \}
\]
be the $\Mod_g$-stabilizer of $\gamma$. Then
\[
 N(\gamma, L) \geq \frac{1}{\# \H} N_\gamma(X,L)
\]
Note that if $\gamma$ is filling then $\H$ must be finite. By \cite{Kerckhoff83}, there is some metric $Z$ on $\S$ so that $\H$ is a subgroup of the isometry group of $Z$ . The Hurwitz automorphism theorem states that the size of a group of conformal automorphisms of a surface $\S$ with hyperbolic metric $Z$ is bounded above by $84(g-1)$ \cite{Hurwitz}. Therefore, we can bound $N(\gamma, L)$ by counting elements of $\Mod_g \cdot X$.
\begin{equation}
\label{eq:FillingComp}
 \frac{1}{84(g-1)} N_\gamma(X,L) \leq N(\gamma, L) \leq  84(g-1)N_\gamma(X,L)
\end{equation}

Instead of counting orbit points in Teichmuller space with respect to the length of $\gamma$, we wish to count orbit points with respect to Teichmuller distance. Let $d_\T(\cdot, \cdot)$ denote Teichmuller distance on Teichmuller space. Let 
\[
 N_\T(X,R) = \# \{f X \in \Mod_g \cdot X \ | \ d_\T(f X, X) < R\} 
\]
be the points in the $\Mod_g$ orbit of $X$ lying in a ball of Teichmuller radius $R$. By \cite[Theorem 1.2]{ABEM12}, 
\[
 N_\T(X,R) \sim c_{X,g}e^{(6g-6)R}
\]
where $c_{X,g}$ is a constant depending on $X$ and the genus $g$ of the surface. (In fact, they give an explicit value for this constant, but we do not need it for this proof.)

If we can relate $N_\gamma(X,L)$ and $N_\T(X,R)$, then we are done. In other words, we need to find the relationship between $l_{f X}(\gamma)$ and $d_\T(f X, X)$. We do this in the case where 
\[
l_X(\gamma) = \min\{l_{fX}(\gamma) \ | \ f \in \Mod_g\}
\]
That is, we suppose that, on $X$, $\gamma$ is the shortest curve in its $\Mod_g$ orbit. We can assume without loss of generality that this is the case, since we can replace $\gamma$ with any curve in $\Mod_g \cdot \gamma$, and the number $N(\gamma, L)$ will be unchanged.
\begin{claim}
\label{cla:filling}
Let $\gamma$ be a filling, closed geodesic that is shortest in its $\Mod_g$ orbit  on $X$. We relate Teichmuller distance to the length of $\gamma$ as follows:
 \[
  \frac{l_{f X}(\gamma)}{l_X(\gamma)} \asymp e^{ d_\T(f X, X)}
 \]
 where the constants only depend on the $Mod_g$ orbits of $X$ and $\gamma$.
\end{claim}

\begin{proof}
%
%

By, for example, \cite[Theorem 2.1]{LRT12}, 
\[
 \left | \log \sup_{\alpha - \text{s.c.c.}}\frac{l_{fX}(\alpha)}{l_X(\alpha)} - d_\T(X, fX) \right | \leq \log c
\]
where the supremum is taken over all simple closed curves and the constant depends only on $X$. (We write the constant as a logarithm to simplify notation later.)

So we just need to compare $\frac{l_{f X}(\alpha)}{l_X(\alpha)}$ with $\frac{l_{f X}(\gamma)}{l_X(\gamma)}$.  By \cite[Proposition 3.5]{Thurston98}, the ratio of lengths is always maximized by simple closed curves:
  \[
   \frac{l_Y(\gamma)}{l_X(\gamma)} \leq \sup_{\alpha - \text{s.c.c.}} \frac{l_Y(\alpha)}{l_X(\alpha)}
  \]
for any metrics $X, Y$. 
Therefore, 
\[
 \frac{l_{fX}(\gamma)}{l_X(\gamma)} \leq  ce^{d_\T(X, fX)}
\]

Next, we show the other direction: $ \frac{l_{f X}(\gamma)}{l_X(\gamma)} \gtrsim e^{d_\T(X, fX)}$.  Because $\gamma$ is filling, \cite[Lemma 5.1]{Basmajian13} implies that
\[
 l_Y(\gamma) \geq \frac 12 \frac{l_Y(\alpha)}{i(\alpha, \gamma)}
\]
for all simple closed curves $\alpha$, and for any hyperbolic metric $Y$. 

We use a trick to bound $i(\alpha, \gamma)$ in terms of $l_X(\alpha)$. Let $\M\L$ and $\P \M\L$ be the spaces of measured, and projective measured, laminations on $\S$, respectively. The function $f : \M \L \rightarrow \R$ with
\[
 f(\alpha) = \frac{i(\alpha, \gamma)}{l_X(\alpha)} 
\]
has the property that $f(c \cdot \alpha) = f(\alpha)$ for all measured laminations $\alpha$ and constants $c \in \R$. Thus, it gives a continuous function $f : \P\M\L \rightarrow \R$. Since $\P\M\L$ is a compact set, there is a constant $d_X$ depending only on $X$ so that
\[
 i(\alpha, \gamma) \leq d_X \cdot l_X(\alpha)
\]
for all measured laminations, and in particular, for all simple closed curves $\alpha$. Therefore,
\begin{align*}
 \frac{l_{f X}(\gamma)}{l_X(\gamma)} &\geq \frac 12 \frac{l_{f X}(\alpha)}{i(\alpha, \gamma)} \cdot \frac{1}{l_X(\gamma)}\\
 &\geq \frac{1}{2d_X \cdot l_X(\gamma)} \cdot \frac{l_{f X}(\alpha)}{l_X(\alpha)}
\end{align*}
Taking the supremum over all simple closed curves $\alpha$, we get 
\[
  \frac{l_{f X}(\gamma)}{l_X(\gamma)}  \geq \frac 1 c \frac{1}{2d_X \cdot l_X(\gamma)}  e^{d_\T(X, fX)}
\]


Since we chose $\gamma$ so that $l_X(\gamma) \leq l_{fX}(\gamma)$ for all $f \in \Mod_g$, the quantity $l_X(\gamma)$ only depends on $X$ and the orbit $\Mod_g \cdot \gamma$. This implies that
\[
  \frac{l_{f X}(\gamma)}{l_X(\gamma)}  \asymp e^{d_\T(X, fX)}
\]
where we only have multiplicative constants, and these constants depend on $X$ and on the orbit $\Mod_g \cdot \gamma$.


\end{proof}

Given this claim, we can finish the proof of Lemma \ref{lem:filling}. Since $\frac{l_{f X}(\gamma)}{l_X(\gamma)} \asymp e^{ d_\T(f X, X)}$, there is some constant $a$ depending only on $X$ and the orbit $\Mod_g \cdot \gamma$ so that 
\[
l_{f X}(\gamma) < L \implies d_\T(f X, X) < \log (a L)
\] 
and
\[
d_\T(f X, X) < \log(\frac 1 a L) \implies l_{f X}(\gamma) < L
\]
So, restricting our attention to the $\Mod_g$ orbit of $X$, we can compare the number of points in a ball around $X$ with the number of points where the length of $\gamma$ is bounded:
\[
 N_\T(X,\log(\frac 1 a L ) \leq N_\gamma(X,L) \leq N_\T(X,\log(a L))
\]
 Thus, inequality (\ref{eq:FillingComp}) implies
\[
\frac{1}{84(g-1)} N_\T(X,\log(\frac 1 a L)) \leq N(\gamma, L) \leq 84(g-1) N_\T(X,\log(a L))
\]

By \cite{ABEM12}, 
\[
 N_\T(X,\log(a L)) \asymp (a L )^{(6g-6)}
\]
and likewise for when $ d_\T(f X, X) < \log(\frac 1 a L)$. So there are constants depending only on $X$ and the orbit  $\Mod_g \cdot \gamma$ so that
\[
 N(\gamma, L) \asymp L^{6g-6}
\]

\end{proof}



\subsection{Non-filling curves}
The case when $\gamma$ is a non-filling curve is similar, but we have to deal with the fact that the $\Mod_g$ stabilizer of $\gamma$ is infinite. So if $\gamma$ is not filling, we only get an upper bound on $N(\gamma, L)$.
\begin{lem}
\label{lem:nonfilling}
 If $\gamma$ is not filling, then
 \[
  N(\gamma, L) \lesssim L^{6g-6}
 \]
where the constants depend only on the $\Mod_g$ orbits of $X$ and $\gamma$.
\end{lem}

We will retrace the steps of the argument for filling curves, and highlight the differences caused by the fact that $\gamma$ is not filling. 

As we did for filling curves, we want to count metrics in the $\Mod_g$ orbit of $X$ instead of counting curves in the $\Mod_g$ orbit of $\gamma$. The obstruction is the existence of the following infinite families of mapping classes. For $f \in \Mod_g$, let 
\[
 [f]_\gamma = \{ g \in \Mod_g \ | \ g \cdot \gamma = f \cdot \gamma \}
\]
Note that $l_{[f]_\gamma X}(\gamma)$ is well-defined, as $l_{f X}(\gamma) = l_{g X}(\gamma)$ for all $g \in [f]_\gamma$. So,
\[
 N(\gamma, L) = \# \{[f]_\gamma \ |\ l_{[f]_\gamma X}(\gamma) \leq L\}
\]

If $\gamma$ fills subsurface $T \subset \S$, then let 
\[
 [f]_T = \{ g \in \Mod_g \ | \ g|_T = f|_T \}
\]
This must be a subset of $[f]_\gamma$. So, since the $\Mod_g$ stabilizer of $X$ has size at most $84(g-1)$ \cite{Hurwitz}, we have
\[
 N(\gamma, L) \leq \# 84(g-1)\{[f]_T X \ |\ l_{[f]_\gamma X}(\gamma) \leq L\}
\]

Define 
\[
 d(X, [f]_T X) = \min_{g \in [f]_T} d_\T(X, g X)
\]
Just as we did with filling curves, we wish to find a relationship between $l_{[f]_\gamma X}(\gamma)$ and $ d(X, [f]_T X)$.


%


\begin{lem}
Let $\gamma$ be the shortest curve in its $\Mod_g$ orbit, on $X$. Then 
\[
 e^{d(X, [f]_T X)} \lesssim \frac{l_{fX}(\gamma)}{l_X(\gamma)}
\]
where the constants depend only on $X$ and $\gamma$.
\end{lem}

\begin{proof}
 By \cite[Theorem B]{CR07}, if $\mu$ is a short marking on $X$, then 
 \[
  \left | d_\T(X,fX) - \log \max_{\alpha \in \mu} \frac{l_{fX}(\alpha)}{l_X(\alpha)} \right | \leq a
 \]
where $a$ depends only on the metric $X$. 

In \cite{CR07}, they define this marking by first choosing a pants decomposition $\alpha_1, \dots, \alpha_{3g-3}$ of $\S$. Then for each $\alpha_i$ in the pants decomposition, they choose a dual curve $\delta_i$ that intersects $\alpha_i$ minimally and is disjoint from $\alpha_j$ for each $j \neq i$. They choose these curves so that the pants decomposition is as short as possible. However, their proof works for any marking of this form. The constant $a$ depends only on the lengths of $\alpha_1, \dots, \alpha_{3g-3}$ and $\delta_1, \dots, \delta_{3g-3}$. So, if for each $\gamma$ that is shortest in its $\Mod_g$ orbit, we find a marking whose length only depends on $l_X(\gamma)$, the constant $a$ given by \cite{CR07} will depend only on the $\Mod_g$ orbit of $\gamma$ and on $X$.

\begin{claim}
\label{cla:ShortPants}
Let $T$ be a surface with geodesic boundary. Then there is a marking $\alpha_1, \dots, \alpha_n, \delta_1, \dots, \delta_n$ so that
\[
l(\alpha_i), l(\delta_i) \leq c_T l_X(\partial T) + c_X, \forall i
\] 
where $c_T$ depends only on the topology of $T$, $c_X$ depends only on the metric $X$ and $l_X(\partial T)$ is the total boundary length of $T$.
\end{claim}
\begin{proof}
This is essentially proven in the proof of \cite[Theorem 5.2.3]{Buser}. Given a surface with boundary, they show that one can construct an arc connecting boundary components, whose length is at most $\epsilon_T = 2\sinh^{-1}\frac{Area(T)}{l_X(\partial T)}$ where $Area(T)$ is determined by the Euler characteristic of $T$ (by the Gauss-Bonnet theorem). From this, we can use the arguments in \cite{Buser} to deduce that the shortest essential simple closed curve in $T$ must have length at most $l_X(\partial T) + 2\epsilon_T$.

Let $\alpha_1$ be this shortest simple closed curve in $T$. Note that $\sinh^{-1}$ is an increasing function. Let
\[
 \epsilon = \min \{\epsilon_T, 2\sinh^{-1}\left(\frac{Area(T)}{l_X(\alpha_1)}\right)\}
\] 
Then, for any subsurface $R \subseteq T$, $\epsilon \geq \epsilon_R$, where $\epsilon_R$ is defined using $l_X(\partial R)$ rather than $l_X(\partial T)$.

In particular, $l_X(\alpha_1) \leq l_X(\partial T) + 2 \epsilon$. We cut $T$ along $\alpha_1$, and get a new surface $T'$. Either $T'$ is a pair of pants, or we get a new shortest simple closed curve $\alpha_2$. 
So, $\alpha_2$ has length at most $l_X(\partial T) + l_X(\alpha_1)  + 2 \epsilon_{T'} \leq 2 l_X(\partial T) + 4 \epsilon$.

Continuing on in this way, we get a pants decomposition $\alpha_1, \dots, \alpha_n$ of $T$, where $\alpha_i$ has length at most $i l_X(\partial T) + 2i \epsilon$. So for each $i$, $l_X(\alpha_i) \leq c_T l_X(\partial T) + c_X$, where $c_T$ depends only on topology of $T$ and $c_X$ depends only on the metric.

We now want to extend the pants decomposition to a marking. Let $\alpha_i$ be a curve in this pants decomposition. It will be contained in at most two pairs of pants. The argument in \cite{Buser} also implies that the shortest arcs from $\alpha_i$ to the boundary of each of these two pairs of pants again have length at most $L + \epsilon$, where $L$ is the total boundary length of the two pairs of pants. Thus, the shortest curve $\delta_i$ that crosses $\alpha_i$ and no other curve in the pants decomposition has length at most $8L + 4 \epsilon$. As $L$ depends only on the length of the curves in the pants decomposition, which in turn depend linearly on $l(\partial T)$ and $\epsilon$, we get the claim.
\end{proof}

\begin{prop}
\label{rem:ActuallyLength}
If $\gamma$ fills $T \subset \S$, then for any metric $Y$ on $T$,
\begin{itemize}
 \item $l_Y(\partial T) \leq 2l_Y(\gamma)$
 \item The shortest marking on $\S$ whose pants decomposition contains $\partial T$ has each curve of length at most \[c_T l_Y(\gamma) + c_Y\] where $c_T$ depends only on the topology of $T$ and $c_Y$ depends only on the metric $Y$.
\end{itemize}
\end{prop}
\begin{proof}
To see the first statement, cut $T$ along $\gamma$. Then for each geodesic boundary component $\beta$ of $T$, $T \setminus \gamma$ contains a cylindrical component $C_\beta$ where one boundary is $\beta$ and the other boundary component is a concatenation of distinct subarcs of $\gamma$, which is homotopic to $\beta$. Thus, $l_X(\beta)$ is bounded above by the total length of these distinct subarcs. Each subarc of $\gamma$ can lie on the boundary of at most two cylinders. Therefore, the total length of $\partial T$ is bounded above by $2l_X(\gamma)$.

Combined with Claim \ref{cla:ShortPants}, we get the statement about the shortest marking on $\S$ containing $\partial T$. The claim does not directly deal with the curves dual to the boundary of $T$. However, its proof implies that these curves must also have length bounded by $c_T l_X(\gamma) + c_X$.
\end{proof}

We will choose the shortest marking $\mu$ of $\S$ with metric $X$, for which the pants decomposition $\alpha_1, \dots, \alpha_{3g-3}$  contains $\partial T$. 
Given $f \in \Mod_g$, we will find a $g \in [f]_T$ so that
\[
 \max_{\alpha \in \mu} \frac{l_{gX}(\alpha)}{l_X(\alpha)} \lesssim \frac{l_{gX}(\gamma)}{l_X(\gamma)}
\]
where the constants will depend only on $X$ and $\Mod_g \cdot \gamma$.

By Proposition \ref{rem:ActuallyLength}, all of the ratios $\frac{l_X(\alpha)}{l_X(\gamma)}$ are bounded from above and below by constants depending only on $X$ and the $\Mod_g$ orbit of $\gamma$. So what we really show is that there is a $g \in [f]_T$ so that for each $\alpha$ in $\mu$, 
\[l_{gX}(\alpha) \lesssim l_{gX}(\gamma)\]

 We can write
\[
 \mu = \mu_T \cup \mu_c
\]
where $\mu_T$ is the set of curves in $\mu$ contained entirely inside $T$ and $\mu_c$ contains the rest of the curves in the marking. Suppose $\alpha \in \mu_T$. Then because $\gamma$ fills $T$, we know by \cite[Lemma 5.1]{Basmajian13} that $l_{fX}(\alpha) \leq i(\alpha, \gamma) l_{fX}(\gamma)$. 
\begin{rem}
 \label{rem:IntersectionIndependent}
 Because $\mu$ was chosen to be the shortest marking containing $\partial T$, and $\gamma$ is the shortest in its $\Mod_g$ orbit, the numbers $i(\alpha, \gamma)$ depend only on $X$ and the $\Mod_g$ orbit of $\gamma$.
\end{rem}

We just need to deal with the case when $\alpha \in \mu_c$. We first find another marking,
\[
 \mu' = \mu_T \cup \mu_c'
\]
that, in fact, satisfies
\begin{equation}
\label{eq:LengthIntersection}
 l_{f X}(\alpha) \leq [i(\alpha, \gamma)  + c_T]l_{f X}(\gamma) + 2c_{fX}
\end{equation}
for each $\alpha \in \mu'$, where $c_T$ and $c_{gX}$ are the constants from Proposition \ref{rem:ActuallyLength}. This should be thought of as the non-filling analogy of \cite[Lemma 5.1]{Basmajian13}, which we used in the case where $\gamma$ is filling.

To form $\mu'$, we need to find $\mu_c'$. The set $\mu_c'$ will contain $\partial T$ and the shortest marking on $\S \setminus T$ in the metric $f X$. By Proposition \ref{rem:ActuallyLength}, all of these curves will have length bounded above by $c_T l_{f X}(\gamma) + c_{fX}$. Thus, they will satisfy (\ref{eq:LengthIntersection}).

In order for $\mu_T \cup \mu_c'$ to be a marking of $\S$, we just need to add the curves that are dual to each boundary curve of $T$. Let $\alpha$ be a boundary curve of $T$, and let $\delta_\alpha$ be its dual curve in the original marking $\mu$. Homotope $\delta_\alpha$ so that in $T$, it is a geodesic arc $\bar \delta_\alpha$ perpendicular to $\alpha$. If we cut $T$ along $\gamma$, we get simply connected regions, and cylindrical regions that correspond to components of $\partial T$. This cuts $\bar \delta_\alpha$ into exactly $i(\delta_\alpha, \gamma)$ pieces. The pieces inside simply connected regions have length at most $l(\gamma)$. (See the proof of \cite[Lemma 5.1]{Basmajian13} for details.) There are two pieces of $\bar \delta_\alpha$ perpendicular to $\partial T$. These lie in the cylindrical regions of $T \setminus \gamma$. The area of each cylindrical region is at most $Area(T)$. So, for example by \cite[Theorem 4.1.1, Theorem 5.2.3]{Buser}, the total length of these pieces of $\bar \delta_\alpha$ are at most $c_{fX}$. Thus, 
\[
l_{fX}(\bar \delta_\alpha) \leq i(\delta_\alpha, \gamma) l_{fX}(\gamma) + c_{fX}
\]

Take the pair of pants $\p$ in the new marking on $\S \setminus T$ that has $\alpha$ as a boundary component.  Now find the shortest arc with respect to metric $f X$ that joins the two endpoints of $\bar \delta_\alpha$ inside $\p$. Call it $\bar \delta_\alpha'$. Then 
\[
l_{fX}(\bar \delta_\alpha') \leq c_T l_{f X}(\gamma) + c_{fX}
\]
The concatenation of $\bar \delta_\alpha$ with $\bar \delta_\alpha'$ gives us the new dual curve $\delta_\alpha'$ to add to $\mu_c'$. By the above, its length is bounded by
\[
 l_{fX}(\delta_\alpha') \leq [i(\delta_\alpha, \gamma) +c_T]l_{fX}(\gamma) + 2c_{fX}
\]
So we have found a set of curves $\mu_c'$ so that the marking $\mu' = \mu_T \cup \mu_c'$ satisfies (\ref{eq:LengthIntersection}).

Because $\Mod_g$ acts cocompactly on the marking graph of $\S$, there is an element $h \in \Mod_g$ and elementary marking moves $m_1, \dots, m_n$ so that $h \cdot \mu = m_n \dots m_1 \mu'$, and the number of marking moves is uniformly bounded in the genus of $\S$. We need to control how each marking move changes the length of the curves. For any metric $Y$, let $L_i$ be the $Y$- length of the longest curve in $m_i \dots m_1 \cdot \mu'$. Then for each $\alpha \in \mu'$,
\[
 l_{Y}(m_{i+1} m_i \dots m_1 \cdot \alpha) \leq 3L_i
\]
This covers the change in length coming from both twist moves (which adds at most $2L_i$ to the length of a transversal) and switch moves, which don't change lengths. Thus, 
\[
 l_{Y}(m_n \dots n_1 \cdot \alpha) \leq 3^n L_0
\]
where $L_0$ is the length of the longest curve in $\mu'$. Therefore,
\[
 l_{fX}(h \mu) \leq 3^{c_\S}[(i(\delta_\alpha, \gamma) +c_T)l_{fX}(\gamma) + 2c_{fX}]
\]
where $c_\S$ is the maximal number of moves that depends only on the topology of $\S$.

Note that
\[
 l_{fX}(h \cdot \mu) = l_{h^{-1}fX} (\mu)
\]
Thus,
\[
 l_{h^{-1}f X} (\alpha) \lesssim l_{h^{-1}f X}(\gamma)
\]
for each $\alpha \in \mu$.
Note that the constants depend on the numbers $i(\alpha, \gamma)$ for $\alpha \in \mu$. But by Remark \ref{rem:IntersectionIndependent}, these numbers depend only on $\Mod_g \cdot \gamma$. Therefore, the constants in the above formula depend only on $X$ and the orbit $\Mod_g\cdot \gamma$.

Without loss of generality, $h \cdot \mu_T = \mu_T$. Thus, $ h^{-1}f\in [f]_T$. So, $l_{h^{-1}f X}(\gamma) = l_{fX}(\gamma)$.  

Let
\[
 g = h^{-1} f
\]
Then,
\[
 \frac{l_{gX}(\alpha)}{l_X(\alpha)} \lesssim  \frac{l_{gX}(\gamma)}{l_X(\gamma)}
\]
for each $\alpha \in \mu$, and where the constants depend only on the $\Mod_g$ orbits of $X$ and $\gamma$. So, by \cite{CR07},
\[
 e^{d(X, [f]_T X)} \lesssim \frac{l_{fX}(\gamma)}{l_X(\gamma)}
\]
where the constants depend only on $X$ and the orbit $\Mod_g \cdot \gamma$.
\end{proof}

This lemma shows that $l_{fX}(\gamma) \leq L$ implies $d(X, [f]_T X) \leq \log (aL + b)$ for some constants $a$ and $b$ that depend on $X$ and the orbit $\Mod_g \cdot \gamma$. Therefore, 
\[
 N(\gamma, L) \leq N_\T(X,\log (aL + b))
\]
By \cite{ABEM12}, this implies that
\[
 N(\gamma, L) \lesssim (aL + b)^{6g-6}
\]
For $L$ large, however, this grows like $L^{6g-6}$. Thus,
\[
 N(\gamma, L) \lesssim L^{6g-6}
\]
for constants depending only on $X$ and the $\Mod_g$ orbit of $\gamma$. This proves Lemma \ref{lem:nonfilling}.

\begin{cor}
\label{cor:UpperBound}
 We deduce the upper bound for $\# \G^c(L,K)$:
 \[
  \# \G^c(L,K) \lesssim L^{6g-6}
 \]
where the constants depend on $K$ and the metric $X$.
\end{cor}
\begin{proof}
 We add up $N(\gamma, L)$ over all $\Mod_g$ orbits of curves with at most $K$ self-intersections. So Lemmas \ref{lem:filling} and \ref{lem:nonfilling} give us this corollary.
\end{proof}

The upper bound in Theorem \ref{thm:Appendix} is given in Corollary \ref{cor:UpperBound}. The lower bound follows from the fact that $\G^c(L,0) \subset \G^c(L,K)$ for all $K$. By \cite{Mirzakhani08}, $\#\G^c(L,0) \sim L^{6g-6}$, so there are some constants so that $L^{6g-6} \lesssim \#\G^c(L,K)$. 

\subsection{Dependence of constants on $K$}
\label{sec:Dependence}
Note that Lemma \ref{lem:filling} also implies that $L^{6g-6} \lesssim \#\G^c(L,K)$. 
The constants we get for the lower bound by using Lemma \ref{lem:filling} are larger than those that come from the lower bound on the number of simple cloesd curves. In fact, we can track dependence of these constants on $K$.

The lower bound in Lemma \ref{lem:filling} follows from the fact that if $d(X, fX) \leq \log(\frac{1}{cl_X(\gamma)} L)$, then $l_{fX}(\gamma) \leq L$. So in the lower bound on $N(\gamma, L)$, the multiplicative constant is bounded below by $\frac{1}{cL_K}$ where $L_K$ is defined as: 
 \[
  L_K = \max_{ \left \{\begin{subarray}{c}
\Mod_g \cdot \gamma \text{ s.t.} \\ i(\gamma, \gamma) \leq K 
        \end{subarray} \right \}
} \min_{f \in \Mod_g}\{l_X(f \cdot \gamma)\} 
 \]
That is, $L_K$ is the largest possible length of the shortest curve in a $\Mod_g$ orbit of a geodesic with at most $K$ self-intersections. We will show in a later paper that the best upper bound on $L_K$ is proportional to $K$.

So we could get a lower bound on $\#\G^c(L,K)$ by summing $N(\gamma, L)$ over all mapping class group orbits of filling curves with at most $K$ self-intersections. If the  number of such orbits is $N_{fill}(K)$, then the multiplicative constant in this lower bound is proportional to $\frac{1}{L_K} N_{fill}(K)$. (The propotionality constant will depend only on the metric $X$.)

Similarly, the multiplicative constant in Corollary \ref{cor:UpperBound} is proportional to $L_K N(K)$, where $N(K)$ is the number of $\Mod_g$ orbits of all curves with at most $K$ self-intersections. Since we do not have a good, explicit estimate for either $N_{fill}(K)$ or $N(K)$, the bounds in Theorem \ref{thm:Appendix} are not explicit in $K$.

\section{Upper bound for arbitrary intersection number}

In what follows, we give an upper bound on $\#\G^c(L,K)$ that is explicit in both $K$ and $L$. To be precise, we show the upper bound on $\#\G^c(L,K)$ for closed surfaces. An upper bound for surfaces with boundary follows as a corollary.

 \begin{theorem}
 \label{thm:ClosedGeneralSurfaceCount}
Let $X_{-1}$ be a negatively curved metric on a closed genus $g$ surface $\S$. For any $L > 0$ and any $K \geq 1$, we get:
 \[
 \# \G^c(L,K) \leq \min 
 \left \{
  c_X e^{\delta L},
  (c_X L)^{c_g \sqrt K} 
 \right \}
\]
where $c_X$ depends on $X_{-1}$, $c_g$ depends only on $\S$, and $\delta$ is the topological entropy of the geodesic flow on $\S$ with respect to $X$.
\end{theorem}
 Note that by Margulis's theorem \cite{MargulisPhD}, $\#\G^c(L,K) \leq c_X e^{\delta L}$ for some constant $c_X$. We just need to prove that $\#\G^c(L,K) \leq (c_X L)^{c_g \sqrt K}$.
 
A corollary of this theorem is an upper bound for $\# \G^c(L,K)$ on a surface with boundary:
 \begin{cor}
Let $X_{-1}$ be a negatively curved metric on a genus $g$ surface $\S$ with $n$ geodesic boundary components. For any $L > 0$ and any $K \geq 1$, we get:
 \[
 \# \G^c(L.K) \leq \min 
 \left \{
  c_X e^{\delta L},
  (c_X L)^{c_\S \sqrt K}
 \right \}
\]
where $c_X$ depends on $X_{-1}$, $c_\S$ depends only on the topology of $\S$, and $\delta$ is the topological entropy of the geodesic flow on $\S$ with respect to $X_{-1}$.
\end{cor}
\begin{proof}
 Let $X_{-1}$ be a negatively curved metric on a genus $g$ surface $\S$ with $n$ geodesic boundary components. We can double it along its boundary to get a closed surface $\S'$ of  genus $2g + n-1$ and negatively curved metric $X_{-1}'$. $\S$ injects into $\S'$ in a conanonical way, so that $X_{-1}'$ pulls back to the metric $X_{-1}$ on $\S$. Thus, closed geodesics in $(\S', X_{-1}')$ either pull back to closed geodesics or multi-arcs in $(\S, X_{-1})$. 
 
 By Theorem \ref{thm:ClosedGeneralSurfaceCount}, $\#\G^c(L,K)$ is at most $(c_{X'} L)^{c_{g'} \sqrt K}$ on $(\S, X_{-1}')$, where $g' = 2g + n-1$, and where $c_{X'}$ is a constant depending on $X_{-1}'$, and therefore on $X_{-1}$. So, on $\S$ with metric $X_{-1}$,
 \[
  \# \G^c(L,K) \leq (c_{X} L)^{c_{\S} \sqrt K}
 \]
where we set $c_X = c_{X'}$ and $c_\S = c_{g'}$.
 
 Furthermore, by extensions of the theorem of Margulis to surfaces with boundary due, for example, to Guillop\'e \cite{Guillope94}, $\# \G^c(L)$ on $(\S, X_{-1})$ is asymptotically $\frac{e^{\delta L}}{\delta L}$, where $\delta$ is the topological entropy of the geodesic flow on $\S$ with respect to $X_{-1}$. Thus, by adjusting the constant $c_X$, we get that $\G^c(L,K) \leq c_X e^{\delta L}$, as well. This gives us the corollary.
\end{proof}

\section{Reduction to flat surfaces}
Theorem \ref{thm:ClosedGeneralSurfaceCount} follows from counting geodesics in a flat metric. To compare geodesics in flat and negatively curved metrics on $\S$, we need to compare their lengths. We do this in the following lemma.
\begin{lem}
 \label{lem:LengthCompatible}
 Let $X$ be a negatively curved metric on $\S$ and let $X_0$ be a flat metric. Then there is a constant $\lambda$ depending on $X$ and $X_0$ so that for all closed geodesics $\gamma \in \G^c$, 
 \[
  \frac 1 \lambda l_0(\gamma) \leq l(\gamma) \leq \lambda l_0(\gamma)
 \]
where $l(\gamma)$ is the length of the geodesic representative of $\gamma$ in $X$ and $l_0(\gamma)$ is the length of the geodesic representative of $\gamma$ in $X_0$. 
\end{lem}
\begin{proof}
 We will use the set $\C(S)$ of geodesic currents for this proof.  As geodesic currents do not appear anywhere else in this paper, we will briefly describe their properties here, and refer the reader to \cite{Bon} for more details. $\C(\S)$ is the set of Borel, geodesic-flow invariant measures on the unit tangent bundle, $T_1(\S)$. The set of closed geodesics $\G^c$ embeds in $\C(\S)$. There is a well-defined intersection number $i(\cdot, \cdot)$ on pairs of geodesics currents that restricts to the usual geometric intersection number on $\G^c \times \G^c \subset \C(\S) \times \C(\S)$. This intersection number is continuous and bi-linear.
 
 Given the negatively curved metric $X$, we can define the associated Liouville current $\mu \in \C(S)$. This geodesic current has the property that for each closed geodesic $\gamma \in \G^c$, $i(\gamma, \mu) = l(\gamma)$. 
 
 By \cite[Theorem 4]{DLR10}, each flat metric $X_0$ on $\S$ can also be represented by a geodesic current $\mu_0 \in \C(S)$. They show that these geodesic currents behave just like the Liouville currents for negatively curved metrics. For example, for each $\gamma \in \G^c$, they show that $i(\gamma, \mu_0) = l_0(\gamma)$. 
 
 Consider the function
\[
\begin{array}{ccccc}
f  &:  & \C(S)             &\longrightarrow     & \R \\ 
  &  &  \gamma   &\mapsto            & \frac{i(\gamma, \mu)}{i(\gamma, \mu_0)}
\end{array}
\]
This map has the property that $f(c \cdot \gamma) = f(\gamma)$ so it descends to a map $f : \P\C(\S) \rightarrow \R$, where $\P\C(\S)$ is the set of projectivized geodesic currents. By \cite[Corollary 5]{Bon}, $\P\C(\S)$ is compact. As $l(\gamma)$ and $l_0(\gamma)$ are never 0,  $f$ is a continuous, positive function on a compact set. Therefore, there are constants $c_1, c_2 >0$ so that 
\[
 c_1 l_0(\gamma) \leq l(\gamma) \leq  c_2 l_0(\gamma)
\]
for each closed geodesic $\gamma \in \G^c$.
\end{proof}

Let $\G^c_Y(L,K)$ denote the set $\G^c(L,K)$ for a metric $Y$ on $\S$. Let $X_{-1}$ be a negatively curved metric on $\S$ and let $X_0$ be a flat metric. Then Lemma \ref{lem:LengthCompatible} implies that
\[
 \G^c_{X_{-1}}(L,K) \subset  \G^c_{X_0}(\frac 1 \lambda L,K)
\]
for each $L, K > 0$. 
Therefore,
\[
 \# \G^c_{X_{-1}}(L,K) \leq  \# \G^c_{X_0}(\frac 1 \lambda L,K)
\]
We give an upper bound on $\#\G^c_{X_0}(L,K)$ in Theorem \ref{thm:ClosedFlatSurfaceBound}. This upper bound directly implies Theorem \ref{thm:ClosedFlatSurfaceBound}.

 \section{Bounding the number of closed geodesics in a flat metric}
 
 Let $X_0$ be a flat metric on $\S$ with one singular point, which we denote $s$. We wish to count closed geodesics with respect to $X_0$ that pass through $s$. If $\gamma$ is  a closed geodesic with respect to $X_0$ that does not pass through $s$, then it is contained in a flat cylinder. But then $\gamma$ must be simple. We know how to count simple closed geodesics, so counting geodesics through $s$ will allow us to count all geodesics on $X_0$. 
 
 Let $\G_*^c$ denote the set of closed geodesics on $X_0$ that are not contained in any cylinder. Then let
 \[
  \G_*^c(L) = \{ \gamma \in \G_*^c \ | \ l_0(\gamma) \leq L\}
 \]
 and
 \[
  \G_*^c(L,K) = \{ \gamma \in \G_*^c \ | \ l_0(\gamma) \leq L, i(\gamma, \gamma) \leq K\}
 \]
Here, $l_0(\gamma)$ denotes the geodesic length of $\gamma$ on $X_0$ and $i(\gamma, \gamma)$ denotes the least transverse self-intersection number of all closed curves in the free homotopy class of $\gamma$.

\begin{theorem}
\label{thm:ClosedFlatSurfaceBound}
Let $X_0$ be a flat metric with one singular point on $\S$. For any $L > 0$ and any $K \geq 1$, we get:
 \[
 \# \G_*^c(L,K) \leq  (c_* L)^{c_g \sqrt K}
\]
where $c_*$ depends on the geometry of $X_0$, and $c_g$ depends only on the topology of $\S$.
\end{theorem}
The bound on all closed geodesics follows directly from this theorem.
\begin{cor}
 Let $X_0$ be a flat metric on $\S$. For all $L > 0$ and all $K \geq 1$, we get:
 \[
 \# \G^c(L,K) \leq  (c_0 L)^{c_g \sqrt K}
\]
 where $c_0$ depends on $X_0$, and $c_g$ depends only on $\S$.
\end{cor}
\begin{proof}
 We know that
 \[
  \G^c(L,K) \setminus \G_*^c(L,K) = \{\gamma \in \G^c \ | \ \gamma \mbox{ lies in a cylinder}\}
 \]
By \cite{Masur90}, there is some universal constant $c_{cyl}$ so that the number of cylinders that contain a closed geodesic of length at most $L$ is at most $c_{cyl} L^2$. Thus,
\[
 \#\G^c(L,K) \leq \#\G_*^c(L,K) + c_{cyl}L^2
\]
Since $\#\G_*^c(L,K) \leq (c_* L)^{c_g \sqrt K}$, and since there is some $L_0$ so that $\G^c(L_0, K) = \emptyset$ for all $K$, there exists a constant $c_X$ depending only on $X$ so that
\[
 \#\G^c(L,K) \leq (c_X L)^{c_g \sqrt K}
\]

\end{proof}

\section{Strategy of the proof}

Let $\C$ be the set of saddle connections on $X_0$. Because there is just one singular point, each $\sigma \in \C$ is a simple arc $\sigma : s \mapsto s$ that corresponds to a simple closed geodesic $\bar \sigma$ that also passes through $s$. Since no geodesic in $\G_*^c$ lies in a flat cylinder, each $\gamma \in \G_*^c$ can be written as the concatenation of saddle connections:
\[
 \gamma = \sigma_1 \dots \sigma_n \mbox{ with } \sigma_i \in \C, \forall i
\]
This should be thought of as a decomposition of $\gamma$ into simple closed curves. Note that the sequence $\sigma_1, \dots, \sigma_n$ uniquely determines $\gamma$.

This is why we work with the flat metric $X_0$ instead of working directly with a negatively curved metric $X$. If, instead, we were working with a hyperbolic metric, it is much more complicated to find an injective map from closed geodesics to collections of simple closed geodesics. For example, for any $\gamma \in \G^c$, one can find many simple closed curves as subarcs of $\gamma$, but they are not all concatenated at a single point. Rather, they are joined by arcs. Given a collection of simple closed geodesics, there are many ways to join them by arcs to get different simple closed curves. We avoid all of these complications by taking closed geodesics in a flat metric with one singular point.


One approach to counting geodesics in $\G_*^c(L,K)$ is as follows. Suppose we can find a function $N(L,K)$ so that if $\gamma = \sigma_1 \dots \sigma_n \in \G_*^c(L,K)$ for $\sigma_i \in \C, \forall i$, then 
\[
n \leq N(L,K)
\]
If $l_0(\gamma) \leq L$, then $l_0(\sigma_i) \leq L, \forall i$. The number of saddle connections of length at most $L$ grows like $O(L^2)$. Thus, 
\[
\# \G_*^c(L,K) \leq cL^{2N(L,K)}
\]
for some constant $c$. 

The problem with this approach is that even simple closed geodesics of length $L$ can be written as roughly length $L$ sequences of saddle connections. So the best we could do is $N(L,K) \approx L$, giving us a bound of $ \#\G_*^c(L,K) \lesssim L^L$. This is not very good. But we get over this problem by replacing sequences of saddle connections with sequences of simple arcs. In particular, the proof goes as follows.

\begin{itemize}
 \item We first define what we mean by a simple geodesic arc $\delta : s \mapsto s$ (Definition \ref{def:SimpleArc}.) 
 
 \item Let
 \[
  \C_0 = \{ \delta : s \mapsto s \ | \ \delta \mbox{ simple geodesic arc} \}
 \]
 and 
 \[
  \C_0(L) = \{ \delta \in \C_0 \ | \ l_0(\delta) \leq L \}
 \]
We bound the size of $\C_0(L)$:
 \[
  \# \C_0(L) \lesssim L^{c_g}
 \]
where $c_g$ is a constant depending only on the genus of $\S$ (Lemma \ref{lem:SimpleCount}).
 \item Lastly, we find a constant $N(L,K)$ so that if $\gamma = \delta_1 \dots \delta_n \in \G_*^c(L,K)$, with $\delta_i \in \C_0, \forall i$, then
 \[
  n \leq N(L,K)
 \]
In fact, 
\[
 N(L,K) \lesssim \min \{ \sqrt K, L\}
\]
(See Lemma \ref{lem:DeltaLengthCount} for the precise statement.)

The fact that a geodesic of length $L$ can be decomposed into a most $L$ simple geodesic arcs is not so surprising. What is interesting is that the number of simple arcs in a geodesic $\gamma$ is also bounded by $\sqrt{i(\gamma, \gamma)}$. 

\item Our theorem then has the form
\[
 \#\G_*^c(L,K) \leq (c_*L)^{c_gN(L,K)}
\]
where $c_*$ is a constant depending on $X_0$.

\end{itemize}

\section{Seeing self-intersections of $\gamma$}

The flat structure $X_0$ on $\S$ gives us a useful decomposition of $\gamma$ into saddle connections. However, geodesics in $X_0$ are generically not self-transverse. So the number of self-intersections of $\gamma$ is not well-defined. We approximate each $\gamma \in \G_*^c$ with a family of nearby curves $\{\gamma_t\}$ so that
\[
 \# \gamma_t \cap \gamma_t = i(\gamma, \gamma), \forall t
\]
In fact, we want to choose $\gamma_t$ to be a geodesic in some negatively curved metric $X_t$, for each $t$. For this, we need the following proposition.

\begin{prop}
 Given a flat metric $X_0$ on $\S$ with one singular point $s$, there is a sequence of negatively curved metrics $\{X_t\}$ so that $\lim_{t \rightarrow 0} X_t = X_0$.
\end{prop}
\begin{proof}
We start by approximating $X_0$ by a sequence of negatively curved metrics with a cone point at $s$. $X_0$ can be formed by gluing together the sides of some $4g$-gon, $A$. This is because we can cut $\S$ along disjoint saddle connections until we get a flat polygon. 

We want to approximate $A$ by $4g$-gons that have constant curvature $-t^2$ for each $0 < t < T$, for some $T$. For each $t > 0$, let $\h_t$ be the plane with constant curvature $-t^2$. Cut $A$ into triangles $T_1, \dots, T_{4g-2}$. For each $t$, take triangles $T^t_1, \dots, T^t_{4g-2}$ in $\h_t$ with the same side lengths as $T_1, \dots, T_{4g-2}$. The side lengths uniquely determine the triangles up to isometry. Thus, $\lim_{t \rightarrow 0} T^t_i = T_i$, for each $i$. Glue the triangles in $\h_t$ together to get a $4g$-gon $A_t$ in $\h_t$ with the same side lengths as $A$. This ensures that $\lim_{t \rightarrow 0} A_t = A$. 

We can glue together opposite sides of $A_t$ by isometries to get a metric $Y_t$ on $\S$. Then $Y_t$ will have constant curvature $-t^2$ outside of the cone point $s$. This cone point has a cone angle that converges to the cone angle of $X_0$. Because $\lim A_t = A$, we get that $\lim_{t \rightarrow 0} Y_t = X_0$ on all compact sets outside of the singular point $s$.

Now for each $t$, we will cut out a disc $D_{cut}(t)$ about $s$, and glue in a smooth disc $D_{glue}(t)$. Let $D_{cut}(t)$ be a disc of radius $3t$. Take local polar coordinates $(r, \theta)$ on $D_{cut}$ so that $s$ lies at $r=0$. We claim that in local coordinates, the metric looks like
\[
 Y_t = dr^2 +  f_t(r) d \theta^2
\]
where
\[
 f_t(r) = \frac{\alpha}{2 \pi} \frac{1}{t} \sinh(tr)
\]
and $\alpha$ is the cone angle at $s$. By, for example, \cite{Petersen06}[Chapter 2, p.47] the curvature of a metric of this form is $- \frac{f_t''(r)}{f_t(r)}$. So we see that the curvature of this metric is $-t^2$. This metric is singular only at $s$. To compute the angle at $s$, we will compute instead the circumference, $c_t$, of a disc of radius $\epsilon$ about $s$:
\[
 c_t = \int_{\theta = 0}^{2 \pi} \frac{\alpha}{2 \pi} \frac{1}{t} \sinh(t\epsilon) d \theta = \alpha \frac 1 t \sinh (t \epsilon)
\]
This is exactly the circumference of a wedge with angle $\alpha$ and radius $\epsilon$. (For example, the circumference of a circle of radius $\epsilon$ in $\h_t$ is $2 \pi \frac 1 t \sinh (t \epsilon)$.) Therefore, this is the correct metric.

Now we want to take a disc $D_{glue}$ with metric $dr^2 + g_t(r) d \theta$ so that  
\begin{itemize}
 \item for some $0 < t_0 < 2t$, $g_t(r)$ satisfies the boundary conditions
 \[
  g_t(r) = \left \{
  \begin{array}{rll}
   \frac 1t \sinh tr & \mbox{ if } & r \in [0, t_0]\\
   \frac{\alpha}{2 \pi} \frac 1t \sinh tr & \mbox{ if} & r \in [2t,3t]
  \end{array}
  \right .
 \]
 and 
 \item $g_t(r)$ is smooth and convex on $[0, 3t]$.
\end{itemize}

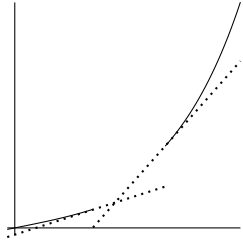
\begin{figure}[H]
 \begin{tikzpicture}
 \draw (-.1,0) -- (3,0); 
 \draw (0,-.1) -- (0,3);
 \draw[domain =-.1:1.04] plot(\x,{.2*sinh(\x)}) ;
  \draw[domain =2:3] plot(\x,{.3*sinh(\x)}) ;
 \draw[thick, dotted,domain = 1.04:3] plot(\x,{1.13 *\x - 1.17});
 \draw[thick, dotted, domain = -.1:2] plot(\x,{.32 *\x - .09});
 \end{tikzpicture}
 \caption{Construction of $h_t(r)$.}
\end{figure}

To construct $g_t(r)$, we just need to join the graph of $\frac 1t \sinh tr$ to the graph of $\frac{\alpha}{2 \pi} \frac 1t \sinh tr$, and get a smooth, convex curve. Draw the tangent line $l_\alpha (r)$ to $f_\alpha(r) = \frac{\alpha}{2 \pi} \frac 1t \sinh tr$ at $t$. Since $f_\alpha(0) = 0$ and $f_\alpha$ is convex, there is some $0 < s_0 < t$ so that $l_\alpha(s_0) = 0$. Draw the tangent line $l_1(r)$ to $f_1(r)$ at $s_0$. Since $f_1$ is convex, and $f_\alpha = \frac{\alpha}{2 \pi} f_1$ for $\alpha > 2 \pi$, there is some time $s_0 < s_1 < t$ so that $l_1(s_1) = l_\alpha(s_1)$. Consider
 \[
  h_t(r) = \left \{
  \begin{array}{rll}
   \frac 1t \sinh tr & \mbox{ if } & r \in [0, s_0]\\
   l_1(r)            & \mbox{ if } & r \in [s_0, s_1]\\
   l_\alpha(r)       & \mbox{ if } & r \in [s_1, t]\\
   \frac{\alpha}{2 \pi} \frac 1t \sinh tr & \mbox{ if} & r \in [t,3t]
  \end{array}
  \right .
 \]
 Then $h_t(r)$ is convex, but not smooth at $s_0, s_1$ or $t$. However, by \cite{Ghomi02}, given any $\delta>0$, there is some function $g_t(r)$ that is smooth and equal to $h_t(r)$ outside of $\delta$ neighborhoods of  $s_0, s_1$ and $t$. In particular, we can find a function $g_t(r)$ and a radius $t_0 < s_0$ so that $g_t(r)$ satisfies the conditions above.

 Because $g_t(r)$ is convex, $D_{glue}$ has negative curvature everywhere. And because of the way that we defined $g_t(r)$, the metric near the boundary of $D_{glue}$ matches up with the metric near the boundary of $D_{cut}$. So we can glue it in to $(\S, Y_t) \setminus D_{cut}$ to get a new negatively curved metric $X_t$ on $\S$ that is smooth at $\partial D_{glue}$.
 
 Next, we see that the angle of $X_t$ at $s$ is $2 \pi$. This is because near $r = 0$, $g_t(r)$ is just like $f_t$ but with $\alpha$ replaced with $2 \pi$. So locally near $s$, $D_{glue}$ looks like a smooth disk of constant curvature $-t^2$.
 
 The last thing we need to check is that the area of $D_{glue}$ goes to zero as $t$ goes to infinity. This will ensure that $\lim_{t \rightarrow 0} X_t = X_0$ on all compact sets outside of $s$. The area of $D_{glue}$ is given by
 \[
 Area_t = \int_0^{3t} \sqrt{ g_t(r) }dr \wedge d \theta
\]
We know that $g_t(r)$ is increasing on $[0,3t]$ because it is a convex function that is increasing at 0. So its maximum value is $\frac{\alpha}{2\pi t} \sinh(3t^2)$. As $\lim_{t \rightarrow 0} \frac 1 t \sinh(3t^2) = 0$, there is some $\epsilon$ small enough so that for all $t < \epsilon$,
\[
 g_t(r) < 1, \forall r \in [0,3t]
\]
Thus, $Area_t < 6 \pi t$ for all $t < \epsilon$. So, $\lim_{t \rightarrow 0}Area_t = 0$. Since $D_{glue}$ is a disc, its radius goes to zero if its area goes to zero. So
\[
 \lim_{t \rightarrow 0}X_t = X_0
\]
on all compact sets outside of $s$, and $X_t$ is a smooth, negatively curved metric for each $t$.

\end{proof}

This proposition allows us to approximate geodesics on $X_0$ by geodesics in nearby negatively curved metrics. 

\begin{lem}
\label{lem:GammaEpsilon}
 For each $\gamma \in \G_*^c$, there is a continuous family of curves $\{\gamma_t\}_{t \in [0,T]}$, with $\gamma_0 = \gamma$ and so that $\gamma_t$ is a geodesic in a negatively curved metric space $X_t$ for each $t$. In particular,
 \[
  i(\gamma, \gamma) = \# \gamma_t \cap \gamma_t, \forall t\in (0,T]
 \]

\end{lem}
\begin{proof}
 Take a sequence of negatively curved metrics $X_t$, where $\lim_{t \rightarrow 0} X_t = X_0$. For each $t$, $\gamma$ is freely homotopic to a closed $X_t$-geodesic $\gamma_t$. So, $\lim_{t \rightarrow 0} \gamma_t = \gamma$, pointwise. Because $\gamma$ has finite length, this limit is, in fact, uniform, and $\{\gamma_t\}_{t \in[0,T]}$ is a continuous family of curves. Geodesics in negatively curved metrics realize self-intersection number, so  $i(\gamma, \gamma) = \# \gamma_t \cap \gamma_t, \forall t$.
\end{proof}

We wish to control how close to $\gamma$ these approximations are. For this we need the following definition.

\begin{defi}
 Two closed curves (or two arcs) $\gamma$ and $\gamma'$ are \textbf{$\epsilon$-homotopic} if there is a homotopy between them that moves each point an $X_0$ distance of at most $\epsilon$. Then we write $\gamma \sim_\epsilon \gamma'$.
\end{defi}
\begin{rem}
  Being $\epsilon$-homotopic is reflexive and symmetric, but not transitive. In fact, if $\gamma_1 \sim_\epsilon \gamma_2$ and $\gamma_2 \sim_\epsilon \gamma_3$ then $\gamma_1 \sim_{2\epsilon} \gamma_3$.
\end{rem}

The problem with approximating a flat geodesic $\gamma$ with an $X_t$-geodesic $\gamma_t$ that realizes its self-intersection number, is that $\gamma_t$ is no longer naturally decomposed into saddle connections. The following lemma gives a way to decompose $\gamma_t$ into approximations of saddle connections.

\begin{figure}[h!] \centering
 \includegraphics{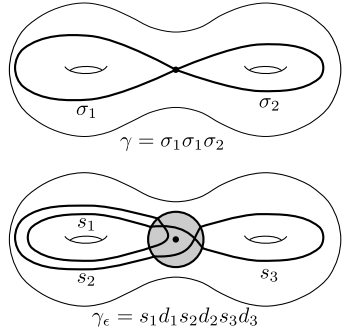}
 \caption[Approximation of a flat geodesic]{How to approximate $\gamma$, while retaining information about saddle connections. The arcs $d_1, d_2$ and $d_3$ lie in the shaded disc $D_{\epsilon'}$.}
 \label{fig:ApproxFlatGeod}
\end{figure}

\begin{lem}
\label{lem:GEpsilonDEpsilon}
 Fix $L$. There is an $\epsilon_L$ depending only on $L$ so that the following holds for all $\epsilon \leq \epsilon' < \epsilon_L$. Let $D_{\epsilon'}$ be an $\epsilon'$-neighborhood of the singular point $s$. For any $\gamma \in \G_*^c(L)$, there is a curve $\gamma_{\epsilon} \sim_{\epsilon} \gamma$ for which $i(\gamma, \gamma) = \# \gamma_{\epsilon} \cap \gamma_{\epsilon}$ and for which we can write 
 \[
  \gamma_{\epsilon} = s_1 \circ d_1 \circ \dots s_n \circ d_n
 \]
where $s_i \subset S \setminus D_{\epsilon'}$ and $d_i \subset D_{\epsilon'}$. 

Furthermore, suppose $ \gamma = \sigma_1 \dots \sigma_n$, with $\sigma_i \in \C, \forall i$. Then for each $i$, 
\[
 s_i \sim_{2\epsilon} \sigma_j \iff  \sigma_j = \sigma_i .
\]
(See Figure \ref{fig:ApproxFlatGeod})
\end{lem}
\begin{proof}
First we choose $\epsilon_L$. For each $\epsilon$, let $D_\epsilon$ be the disc of radius $\epsilon$ about $s$. Because there are finitely many saddle connections of length at most $L$, there is some distance $\epsilon_L'$ so that if $\epsilon < \epsilon_L'$ and if $\sigma \in \C$ with $l_0(\sigma) \leq L$, then $\sigma$ crosses $\partial D_\epsilon$ exactly twice. This means that $\sigma$ will not dip multiple times into $\D_{\epsilon}$ as it travels around $\S$. Let $l_0$ be the length of the shortest closed geodesic on $X_0$. Then we set
\[
\epsilon_L = \min \{ \epsilon_L', \frac{l_0}{8}\} 
\]

Choose $\epsilon$ and $\epsilon'$ so that $\epsilon \leq \epsilon' < \epsilon_L$, and take the disc $D_{\epsilon'}$.

Take a continuous family $\{\gamma_t\}$, for $t \in [0, T]$, where $\gamma_0 = \gamma$ and where, for each $t$, $\gamma_t$ is a geodesic in negatively curved metric $X_t$ on $\S$ (from Lemma \ref{lem:GammaEpsilon}.) There is some $t_0$ depending on $\epsilon$ so that for all $t \leq t_0$,  $\gamma_t \sim_{\epsilon} \gamma$, and the homotopy on $[0, t_0]$ is transverse to $\partial D_{\epsilon'}$. Thus, the number of intersections of $\gamma_t$ with $\partial D_{\epsilon'}$ remains constant for all $t \in [0, t_0]$. Define 
\[
 \gamma_{\epsilon} = \gamma_{t_0}
\]

Write $\gamma = \sigma_1 \dots \sigma_n$, for $\sigma_i \in \C, \forall i$. For each $\sigma_i$, the homotopy $\{\gamma_t\}_{t \in [0,t_0]}$ gives a correspondence between $\sigma_i$ and a subarc $(\sigma_i)_{\epsilon}$ of $\gamma_{\epsilon}$. Since the homotopy moves each point of $\gamma$ by at most $\epsilon < \epsilon'$, the endpoints of $(\sigma_i)_{\epsilon}$ lie inside $D_{\epsilon'}$. Because $\epsilon < \epsilon_L$, each saddle connection crosses $\partial D_{\epsilon'}$ exactly twice. Because the homotopy $\{\gamma_t\}$ is transverse to $\partial D_{\epsilon'}$, the arc $(\sigma_i)_{\epsilon}$ also crosses $\partial D_{\epsilon'}$ exactly two times. Let $s_i$ denote the part of $(\sigma_i)_{\epsilon}$ outside of $D_{\epsilon'}$. We get that 
\[
 \gamma_\epsilon = s_1 d_1 \dots s_n d_n
\]
where $s_1, \dots, s_n$ are the arcs defined above, and $d_i$ connects $s_{i-1}$ to $s_i$. Because $\gamma_\epsilon$ crosses $\partial D_{\epsilon'}$ only at the endpoints of $s_1, \dots, s_n$, each $d_i$ must be contained inside $D_{\epsilon'}$.

Now we need to show that $s_i \sim_{2 \epsilon} \sigma_j$ if and only if $\sigma_j = \sigma_i$. Because $\sigma_i \sim_{\epsilon} (\sigma_i)_{\epsilon}$ and because $s_i \sim_\epsilon (\sigma_i)_\epsilon$, we have that $\sigma_i \sim_{2\epsilon} s_i$. 

Suppose $\sigma_j \sim_{2\epsilon} s_i$ for some $j$. Then $\sigma_j \sim_{4 \epsilon} \sigma_i$. If $\sigma_i \neq \sigma_j$, the $4 \epsilon$ homotopy between them sends some endpoint of $\sigma_i$ to an endpoint of $\sigma_j$ along a non-trivial loop based at $s$. This loop can be tightened to a closed geodesic, whose length must be at least $l_0$. By assumption, $4 \epsilon < \frac{l_0}{2}$. Thus, two saddle connections are $4\epsilon$-homotopic if and only if they are equal. Therefore, $\sigma_j \sim_{2\epsilon} s_i$ if and only if $\sigma_i = \sigma_j$. 
\end{proof}

\section{Intersection number for arcs}
Take a geodesic arc $\delta: s \mapsto s$. We want to define a geodesic self-intersection number $i(\delta,\delta)$ that is intrinsic to $\delta$. This intersection number should have the following property. 

\begin{figure}[h!]
 \centering
 \includegraphics{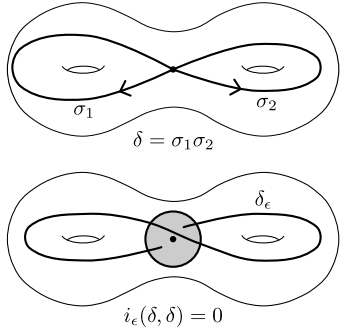}
 \caption{The arc $\delta$ is simple, even though the closed curve given by $\sigma_1 \sigma_2$ has one self-intersection.}
 \label{fig:ArcIntNum}
\end{figure}

Suppose  $\gamma$ is a geodesic in $X_0$, and take some curve $\gamma_\epsilon \sim_\epsilon \gamma$. Suppose our arc $\delta$ happens to be a subarc of $\gamma$. The homotopy from $\gamma$ to $\gamma_\epsilon$ gives a correspondence between $\delta$ and some subarc $\delta_\epsilon$ of $\gamma_\epsilon$. As long as $\epsilon$ is small enough, we want
\[
 i(\delta, \delta) \leq \# \delta_\epsilon \cap \delta_\epsilon
\]
where the left hand side is the intrinsic self-intersection number defined below, and the right hand side is the number of intersections we observe in $\delta_\epsilon$. This is formalized in the following definition:

\begin{defi}
\label{def:SimpleArc}
 Let $\delta_1, \delta_2 : s \mapsto s$ be two geodesic arcs in $X_0$. For each $\epsilon > 0$, let
 \[
  i_\epsilon(\delta_1, \delta_2) = \inf \{ \# (\delta_1)_\epsilon \cap (\delta_2)_\epsilon \ | \ \delta_i \sim_\epsilon (\delta_i)_\epsilon, i = 1,2 \}
 \]
 and let
 \[
  i(\delta_1, \delta_2) = \lim_{\epsilon \rightarrow 0} i_\epsilon(\delta_1, \delta_2)
 \]
 Note that when $\delta_1 = \delta_2$, we require $(\delta_1)_\epsilon = (\delta_2)_\epsilon$, and we count the number of transverse self-intersections of $(\delta_1)_\epsilon$. Thus, a \textbf{simple geodesic arc} is one whose self-intersection number is zero in this sense. 
(See Figure \ref{fig:ArcIntNum}.)
\end{defi}

\begin{rem}
\label{rem:IncFunc}
 The limit $\lim_{\epsilon \rightarrow 0} i_\epsilon(\delta_1, \delta_2)$ exists: If we have some geodesic arc $\delta$ and two values $\epsilon < \epsilon'$, then $\delta_{\epsilon} \sim_{\epsilon} \delta$ implies $\delta_{\epsilon} \sim_{\epsilon'} \delta$, too. For this reason, $i_\epsilon (\delta_1, \delta_2)$ is an increasing function of $\epsilon$. It is bounded above because the arcs have finite length. Thus the limit $\lim_{\epsilon \rightarrow 0} i_\epsilon(\delta_1, \delta_2)$ must exist.
\end{rem}

\section{Counting simple arcs}

Let $\C_0$ be the set of simple geodesic arcs:
\[
 \C_0 = \{ \delta: s \mapsto s \mbox{ geodesic} \ | \ i(\delta, \delta) = 0 \mbox{ as an arc}\}
\]
and consider those simple arcs of length less than $L$:
\[
 \C_0(L) = \{ \delta \in \C_0 \ | \ l_0(\delta) \leq L\}
\]

\begin{lem}
\label{lem:SimpleCount}
 Fix an $L>0$. Then we get the following upper bound on the size of $\C_0(L)$:
 \[
  \# \C_0(L) \leq    (c_0 L)^{c_g}
 \]
where $c_0$ is a constant depending only on the geometry of $X_0$ and $c_g$ is a constant depending only on the surface $\S$.
\end{lem}
\begin{proof}

 We first fix an $\epsilon'$ for the proof of this lemma (and for all the claims used to prove it). There is some $\mu_L>0$ depending only on $L$ so that $\forall \epsilon' < \mu_L$, $\forall \sigma_1, \sigma_2 \in \C$ with $l_0(\sigma_i) < L, i = 1,2$, the two saddle connections $\sigma_1$ and $\sigma_2$ do not intersect on $\partial D_{\epsilon'}$. Such a $\mu_L$ exists because the set $\{\sigma \in \C \ | \ l_0(\sigma) <L\}$ is finite, so the set of intersection points between pairs of saddle connections in this set is also finite. So there is some $\mu_L > 0$ so that $D_{2 \mu_L}$ contains none of these intersection points. Let $\epsilon_L$ be the constant from Lemma \ref{lem:GEpsilonDEpsilon}. Now choose any
 \[
  \epsilon' < \min\{\epsilon_L, \mu_L\}
 \]

 The proof of Lemma \ref{lem:SimpleCount} goes as follows. 
 \begin{itemize}
  \item We fix some set $\Sigma = \{\sigma_1', \dots, \sigma_m'\}$ of distinct saddle connections, and consider those $\delta \in \C_0(L)$ composed only of saddle connections in $\Sigma$. That is, we set
  \[
   \C_0(L,\Sigma) =  \{ \delta = \sigma_1 \dots \sigma_n \in \C_0(L) \ | \ \forall i, \sigma_i \in \Sigma \}
  \]
 
 Then we bound the size of this set (Claim \ref{cla:FixedSaddleConn}.) We use techniques that are similar to those found in \cite{BS85}. Roughly, a geodesic in $\C_0(L,\Sigma)$ is given by weights on the arcs in $\Sigma$, together with data that give the order in which the arcs are traversed. 
 
 The bound on length gives restrictions on which weights are possible. The fact that arcs in $C_0(L,\Sigma)$ are simple restricts the order in which the saddle connections can be traversed. This allows us to bound the size of $C_0(L,\Sigma)$.
 
 \item We want to understand the sets $\Sigma = \{\sigma_1', \dots, \sigma_m'\}$ of distinct saddle connections that can form in arcs in $\C_0(L)$.
 
 We show that $\# \Sigma \leq b_g$ for some constant $b_g$ depending only on the topology of $\S$  (Claim \ref{cla:NumSaddleConn}). This follows from the fact that if $\sigma_1', \dots, \sigma_m'$ all appear in a simple arc, then as arcs, these saddle connections can all be realized disjointly. So the number of such arcs depends only on the genus $g$ of $\S$.
 
 
 \item Lastly, we bound the number of different sets $\Sigma$ of distinct saddle connections that form simple arcs in $\C_0(L)$ (Claim \ref{cla:HowManySigma}). This follows from combining the bound on $\# \Sigma$ with the fact that there are at most $ O(L^2)$ saddle connections of length at most $L$ on $X_0$.
 
 \item We sum our bounds on $\#\C_0(L,\Sigma)$ over all possible sets $\Sigma$ to get an upper bound on $\#\C_0(L)$.
 \end{itemize}

 \begin{claim}
 \label{cla:FixedSaddleConn}
  Let $\Sigma = \{\sigma'_1, \dots, \sigma'_m\}$ be a set of $m$ distinct saddle connections with $l_0(\sigma_i')\leq L, \forall i$. 
  Let
  \[
   \C_0(L,\Sigma) = \{ \delta = \sigma_1 \dots \sigma_n \in \C_0(L) \ | \ \forall i, \sigma_i \in \Sigma\}
  \]
 be the set of $\delta \in \C_0(L)$ composed of the saddle connections in $\Sigma$. Then
 \[
  \#\C_0(L, \Sigma) \leq 16 \left (\frac{L}{l_0} \right )^{m^2+4}
 \]
 where $l_0$ is the length of the shortest closed geodesic on $X_0$.
 \end{claim}
 \begin{proof}
 This argument is inspired by techniques from the proof of a theorem of Birman and Series \cite{BS85}.
 

 
Consider the points $x_1, \dots, x_{2m}$ where the saddle connections in $\Sigma$ intersect $\partial D_{\epsilon'}$.  We have fixed an $\epsilon' < \min\{\epsilon_L, \mu_L\}$ at the start of the proof of Lemma \ref{lem:SimpleCount}. Because $l_0(\sigma_i') \leq L$ for each $i$, our choice of $\epsilon'$ guarantees that $x_1, \dots, x_{2m}$ are all distinct. Let $I_r(x_i)$ be the ball of radius $r$ about $x_i$ in $\partial D_{\epsilon'}$. Choose $r$ small enough so that $I_r(x_i)$ and $I_r(x_j)$ are disjoint for each $i \neq j$.  From now on, let
 \[
  I_i = I_r(x_i)
 \]
 
 Suppose $\delta \in \C_0(L, \Sigma)$. Write $\delta = \sigma_1 \dots \sigma_n$, for $\sigma_i \in \Sigma, \forall i$. Let $\epsilon = \min\{ \frac r2, \epsilon'\}$. The proof of Lemma \ref{lem:GEpsilonDEpsilon} never used the fact that $\gamma$ was closed. By assumption, $\epsilon \leq \epsilon' < \epsilon_L$, where $\epsilon_L$ is the number from Lemma \ref{lem:GEpsilonDEpsilon}. So that lemma implies that there is an arc $\delta_\epsilon \sim_{\epsilon} \delta$ so that $\# \delta_\epsilon \cap \delta_\epsilon = 0$ and that we can write as
 \[
  \delta_\epsilon = s_1 d_1 \dots d_{n-1} s_n
 \]
 where $s_i$ lies outside the disc $D_{\epsilon'}$, $s_i \sim_{2\epsilon} \sigma_i$, and $d_i$ lies inside $D_{\epsilon'}$, for each $i$. Note that we cut off the small subarcs at the ends of $\delta_\epsilon$ that lie inside $D_{\epsilon'}$. (See Figure \ref{fig:ArcIntervals}, but ignore the caption for now.)
 
 Suppose $\sigma_i$ has endpoints $x_{j_i}$ and $x_{k_i}$.  Since $s_i \sim_{2\epsilon} \sigma_i$, and $\epsilon <\frac r2$, the endpoints of $s_i$ lie in $I_{j_i}$ and $I_{k_i}$. Thus, the arcs $s_1, \dots, s_n$ connect the intervals $I_1, \dots, I_m$ outside $D_{\epsilon'}$, and the arcs $d_1, \dots, d_{n-1}$ connect these intervals inside $D_{\epsilon'}$ (Figure \ref{fig:ArcIntervals}).
 
 \begin{figure}[h!] \centering
  \includegraphics{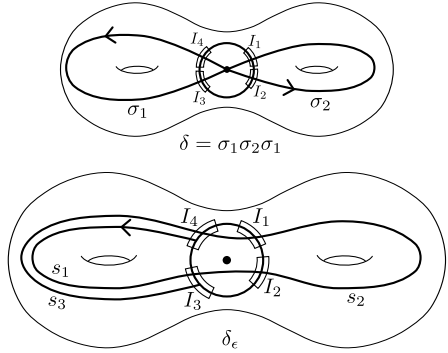}
  \caption[Example of the construction in Claim \ref{cla:FixedSaddleConn}]{In this example, $n_{14}=1, n_{23}=1$ and the other $n_{ij} =0$. If we count from the top corner of $I_1$, $t_0 = 5$ and $t_1 = 3$. Lastly, $i_0 = 4$ and $i_1 = 3$.}
  \label{fig:ArcIntervals}
 \end{figure}

 Now we are ready to give combinatorial data that encodes how many times each $\sigma_i'$ appears in $\delta$, as well as the order of the saddle connections inside $\delta$. Let $n_{ij}$ be the number of arcs that connect $I_i$ to $I_j$ inside $D_{\epsilon'}$. Number the intersection points of $\delta_\epsilon$ with $\partial D_{\epsilon'}$ clockwise from some fixed endpoint of $I_1$. Let $t_0$ and $t_1$ be the number of the start- and endpoints of $\delta_\epsilon$, respectively. Let $I_{i_0}$ and $I_{i_1}$ be the intervals that contain the start- and endpoints of $\delta_\epsilon$, respectively. Then $\delta$ has data $D(\delta) = \big \{ \{n_{ij}\}, t_0, t_1, i_0, i_1 \big \}$. (See Figure \ref{fig:ArcIntervals}.)

 We will show that only $\delta$ can have data $D(\delta)$. First of all, the data determines the number of times each saddle connection $\sigma_i' \in \Sigma$ appears in $\delta$. For each $i$, let
  \[
  n_i =  \sum_j n_{ij}
 \]
Suppose $\sigma_k'$ has an endpoint on $I_{j_k}$. If $j_k \neq i_0, i_1$, then $\sigma_k'$ appears $n_{j_k}$ times. Otherwise, it appears either $n_{j_k}+1$ or $n_{j_k}+2$ times, depending on whether just one of $i_0$ and $i_1$ is $j_k$, or if $i_0 = i_1 = j_k$, respectively.

If $\delta'$ has the same data as $\delta$ then we have shown that the saddle connections in $\Sigma$ appear in $\delta$ and $\delta'$ with the same multiplicity. We need to show that the saddle connections also appear in the same order. This will imply $\delta = \delta'$. 

Take an arc $\delta'_\epsilon \sim_\epsilon \delta'$, that we can write $\delta'_\epsilon = s_1' d_1' \dots d_{n-1}' s_n'$, where $s_i'$ lies outside $D_{\epsilon'}$ and $d_i'$ lies inside $D_{\epsilon'}$. Suppose $\delta'_\epsilon$ intersects $\partial D_{\epsilon'}$ at points $y_0, \dots, y_{2n}$. Suppose the indices on these points correspond to their order around $D_{\epsilon'}$. We will show that we can recover the order in which $y_0, \dots, y_{2n}$ appear in $\delta'_\epsilon$ just from the data.

The arc $\delta'_\epsilon$ gives a pairing of the set of points $\{y_0, \dots, y_{2n}\} \setminus \{y_{t_0}, y_{t_1}\}$ by arcs inside $D_{\epsilon'}$ and a pairing of the set of points $\{y_0, \dots, y_{2n}\}$ by arcs outside $D_{\epsilon'}$. We will actually show that we recover both of these pairings. This will give us the order in which $y_0, \dots, y_{2n}$ appear in $\delta'_\epsilon$. 

Since $\delta'$ has the same data as $\delta$, it also has $n_{ij}$ of the arcs in the set $\{d_1', \dots, d_{n-1}'\}$ connecting points on $I_i$ to points on $I_j$. Because $D_{\epsilon'}$ is a disc, there is only one way to pair the points by \textit{disjoint} arcs inside $D_{\epsilon'}$ so that $n_{ij}$ points on $I_i$ are joined to points on $I_j$. Therefore, the data determines the pairing of points inside $D_{\epsilon'}$.

Now we turn to the pairing of points by arcs outside $D_{\epsilon'}$. Suppose a saddle connection $\sigma_i' \in \Sigma$ joins interval $I_{j_i}$ to interval $I_{k_i}$. Then the points in $\{y_0, \dots, y_{2n}\}$ that lie on $I_{j_i}$ can only be paired to those points that lie on $I_{k_i}$. These points are paired by a set of disjoint arcs outside $D_{\epsilon'}$ that are $\epsilon$-homotopic to $\sigma_i'$. We claim that only one pairing by disjoint arcs is possible. The proof of this is a bit technical, but roughly speaking all we do is lift everything to the universal cover to reduce this to a problem of connecting points on the boundary of a simply connected domain. (See Figure \ref{fig:SigmaiNbhd}.)

Let $N_\epsilon(\sigma_i')$ be an $\epsilon$-neighborhood of $\sigma_i'$. If $s'_j \sim_{2 \epsilon} \sigma_i'$, then $s'_j \in N_\epsilon(\sigma_i')$. (This follows from the construction of $s'_1, \dots, s'_n$ in Lemma \ref{lem:GEpsilonDEpsilon}.) Lift $N_\epsilon(\sigma_i')$ to a region $\tilde N_\epsilon(\sigma_i')$ in the universal cover. This is an $\epsilon$-neighborhood of some lift $\tilde \sigma_i'$ of $\sigma_i'$. There are two lifts $(\tilde D_\epsilon)_1$ and $(\tilde D_{\epsilon})_2$ at either end of $\tilde N_\epsilon(\sigma_i')$. Because $\epsilon < \epsilon_L$, $(\tilde D_\epsilon)_1$ and $(\tilde D_{\epsilon})_2$  are disjoint. In fact, $\tilde N_{\epsilon}(\sigma_i')$ is composed of $(\tilde D_\epsilon)_1$, $(\tilde D_{\epsilon})_2$, and a simply connected region $\tilde R$ between them (see Figure \ref{fig:SigmaiNbhd}). As $s'_1, \dots, s'_n$ lie outside of $D_\epsilon$ but inside $N_\epsilon(\sigma_i')$, their lifts $\tilde s'_1, \dots, \tilde s'_n$ lie in $\tilde R$ and have endpoints on the boundaries of $(\tilde D_\epsilon)_1$ and $(\tilde D_{\epsilon})_2$. Because $\tilde R$ is simply connected, is just one way to join the endpoints of $\tilde s'_1, \dots, \tilde s'_n$ lying on  $\partial (\tilde D_\epsilon)_1$ and $\partial (\tilde D_{\epsilon})_2$.

 
 \begin{figure}[h!] \centering
  \includegraphics{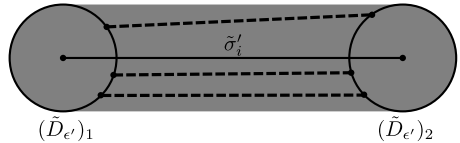}
  \caption[The unique way to connect the intervals with non-intersecting lines outside of $D_{\epsilon'}$]{The $\epsilon$-neighborhood of $\tilde \sigma_i'$ is shaded. There is only one way to join the points on $\partial (\tilde D_{\epsilon'})_1$ to the points on $\partial  (\tilde D_{\epsilon'})_2$ inside the shaded region.}
  \label{fig:SigmaiNbhd}
 \end{figure}
 
 Therefore, the pairings on the sets $\{y_0, \dots, y_{2n}\} \setminus \{y_{t_0}, y_{t_1}\}$ and $\{y_0, \dots, y_{2n}\}$ by arcs inside and outside of $D_{\epsilon'}$, respectively, is determined by the data $D(\delta)$. 
 
 So we have shown that, given any arc $\delta'$ with data $D(\delta)$, and any $\delta'_\epsilon \sim_\epsilon \delta'$, we know how many crossing points $y_0', \dots, y_{2n}'$ $\delta_\epsilon'$ must have, and which intervals these crossing points are contained in. If we choose an arbitrary cyclic order of these points, this also determines the order in which these crossing points are traversed by $\delta'_\epsilon$.
 
 Suppose we index the points $y_0, \dots, y_{2n}$ in the order in which they appear in $\delta_{\epsilon}$. This tells us the order in which the saddle connections appear in $\delta$. Each point $y_j$ lies in some interval $I_i$ and each interval contains the endpoint of a unique saddle connection in $\Sigma$. So each point $y_j$ corresponds to the endpoint of some saddle connection in $\Sigma$. In fact, for each $k = 1, \dots, n$, each pair $(y_{2j}, y_{2j+1})$ corresponds to a pair of intervals that contain the endpoints of the same saddle connection, $\sigma_{i_j}$. Thus, the pair $(y_{2j}, y_{2j+1})$ correspond to $\sigma_{i_j}$ with a specific orientation. So the ordered sequence of points $y_0, \dots, y_{2n}$ tell us that $\delta = \sigma_1 \dots \sigma_n$, where the saddle connections appear with the appropriate orientation. Therefore, the data $D(\delta)$ uniquely determines the simple arc $\delta$.

 Each saddle connection on $X_0$ has length at least $l_0$. As $l_0(\delta) \leq L$ for each $\delta = \sigma_1 \dots \sigma_n \in \C_0(L, \Sigma)$, the number $n$ of saddle connections must be at most $\frac{L}{l_0}$. Given data $D(\delta) = \big \{\{n_{ij}\}, t_0, t_1, i_0, i_1 \big \}$ for some $\delta \in \C_0(L, \Sigma)$, this implies that 
 \[
  \sum_{i,j = 1}^m n_{ij} \leq \frac{L}{l_0}
 \]
 where $m$ is the size of the set $\Sigma$. This sum has at most $m^2$ terms. Thus the number of sets $\{n_{ij}\}$ that satisfy this inequality is at most $(\frac{L}{l_0})^{m^2}$. Furthermore, $t_0, t_1, i_0$ and $i_1$ are integers between 1 and $\frac{2L}{l_0}$, so there are at most $(\frac{2L}{l_0})^4$ choices for them. Thus, the number of possible sets of data given sigma is at most
 \[
 2^4  \left (\frac{L}{l_0} \right )^{m^2+4}
 \]
\end{proof}

We now give some restrictions on the sets of distinct saddle connections that can be used to form simple arcs.

 \begin{claim}
\label{cla:NumSaddleConn}
  Let $\Sigma = \{\sigma'_1, \dots, \sigma'_m\}$ be the set of distinct saddle connections that appear in some $\delta \in \C_0$. Then the size of the set is bounded by 
  \[
   m \leq b_g
  \]
 where $b_g = 2g-1$.
 \end{claim}
 \begin{proof}
 Because $X_0$ has just one singular point, if two different saddle connections intersect outside of $s$, then they intersect transversally. Suppose $\delta $ is simple, and can be written as $\delta = \sigma_1 \dots \sigma_n$, for $\sigma_i \in \C, \forall i$. Then $\sigma_1, \dots, \sigma_n$ cannot intersect outside of $s$. If they did, then for the constant $\epsilon' $ defined at the begining of the proof of Lemma \ref{lem:SimpleCount}, take the disk $D_{\epsilon'}$ about $s$. If $\sigma_1$ and $\sigma_2$ intersect, then $\sigma_1, \sigma_2$ and $\partial D_{\epsilon'}$ form a loop. Because $\sigma_1$ and $\sigma_2$ are geodesics, this loop cannot be homotopic to a point. For any $\epsilon < \epsilon'$, if $\delta_\epsilon \sim_\epsilon \delta$, then the endpoints of $\delta_\epsilon$ will stay inside $D_\epsilon$, so this loop will persist. But then, $\delta_\epsilon$ will not be simple for any $\epsilon$, meaning $\delta$ is not simple.
 
 Suppose $\Sigma = \{\sigma_1', \dots, \sigma_m'\}$ is the set of distinct saddle connections that appear in $\delta$. Without loss of generality, it is a maximal set of disjoint saddle connections (otherwise, we can add saddle connections to $\Sigma$ until this is true.) This gives a cell decomposition of $\S$ with one vertex and $m$ edges. An Euler characteristic argument implies that $m \leq 2g-1$.

 \end{proof}

 \begin{claim}
 \label{cla:HowManySigma}
  The number of sets $\Sigma = \{\sigma_1', \dots, \sigma_m'\}$ that can occur as sets of distinct saddle connections in arcs $\delta \in \C_0(L)$ is at most
  \[
   (b_0 L)^{2b_g}
  \]
 where $b_g$ is a constant depending only on $\S$, and $b_0$ is a constant depending on $X_0$.
 \end{claim}
 \begin{proof}
If a saddle connection $\sigma$ occurs in some $\delta \in \C_0(L)$, then $l_0(\sigma) \leq L$. The number of saddle connections on $X_0$ of length at most $L$ is bounded above by $b_0 L^2$, where $b_0$ is a constant depending only on $X_0$ \cite{Masur90}. For any set $\Sigma$ of distinct saddle connections occurring in some $\delta \in \C_0(L)$, $\#\Sigma \leq b_g$. The number of ways to choose at most $b_g$ elements from a set of size $b_0L^2$ is at most ${2b_0L^2 \choose b_g}$. To see this, suppose we have a list containing two copies of each element of $\C_0(L)$. We choose $b_g$ of the elements in this list. But, we keep an element we choose if and only it is the first time that element appears in our list. This gives us at most $b_g$ elements of $\C_0(L)$. Note that ${2b_0L^2 \choose b_g}$ is bounded above by
\[
 (2 b_0L^2)^{b_g}
\]
 \end{proof}

 Combining Claims \ref{cla:FixedSaddleConn} and \ref{cla:NumSaddleConn}, we get that 
 \[
  \#\C_0(L, \Sigma) \leq 16 \left (\frac{L}{l_0} \right )^{d_g}
 \]
for each set $\Sigma$ of distinct saddle connections that occur in some $\delta \in \C_0(L)$, and where $d_g = b_g^2+4$ is a constant depending only on $\S$. By Claim \ref{cla:HowManySigma}, there are at most $(b_0 L)^{2b_g}$ choices for $\Sigma$, so summing $\#\C_0(L,\Sigma)$ over all sets $\Sigma$ we get that 
 \[
  \# \C_0(L) \leq c_0 L^{c_g}
 \]
 where $c_0 =16 b_0^{2b_g}/l_0^{d_g}$ and $c_g = b_g^2 + 2 b_g + 4$.
\end{proof}

\section{Bounding the number of simple arcs in a closed geodesic}

Suppose $\gamma = \sigma_1 \dots \sigma_n \in \G_*^c$, for $\sigma_i \in \C$ for each $i$. Then we can take the smallest partition of $\sigma_1, \dots, \sigma_n$ into simple arcs:
\[
 \gamma = \delta_1 \dots \delta_m
\]
for $\delta_1 = \sigma_1 \dots \sigma_{n_1}$, $\delta_2 = \sigma_{n_1+1} \dots \sigma_{n_2}$, and so on, with $\delta_i \in \C_0$ for each $i$. In particular, if we concatenate $\delta_i$ and $\delta_{i+1}$, the arc $\delta_i \delta_{i+1}$ is not simple. 

If $\gamma \in \G_*^c(L)$, then each of the arcs $\delta_1, \dots, \delta_m$ has length at most $L$. We know how to count simple arcs of length at most $L$, so we just need to bound the length $m$ of the sequence of these arcs in terms of $l_0(\gamma)$ and $i(\gamma, \gamma)$.

\begin{lem}
\label{lem:DeltaLengthCount}
 Let $\gamma \in \G_*^c(L,K)$, with $K \geq 1$. Suppose
 \[
  \gamma = \delta_1 \dots \delta_m
 \]
is the shortest way to write $\gamma$ as a concatenation of arcs $\delta_i \in \C_0$. Then
\[
 m \leq \min\{\frac{L}{l_0},c \sqrt K\}
\]
where $l_0$ is the length of the shortest closed geodesic on $X_0$ and $c$ is a constant depending only on the topology of $\S$.
\end{lem}
\begin{proof}
 Bounding $m$ in terms of $l_0(\gamma)$ is relatively simple. The difficulty lies in bounding $m$ in terms of $i(\gamma, \gamma)$, which we do first.
 
 Suppose $\gamma \in \G_*^c(L,K)$ with $\gamma = \delta_1 \dots \delta_m$ and $\delta_i \in \C_0$, $\forall i$. Suppose this is the shortest way to represent $\gamma$ as a concatenation of simple arcs. Then, as previously mentioned, the arc $e_i = \delta_i \delta_{i+1}$ is not simple. In fact, let
 \[
  \C_1 = \{ e = d d' \ | \ d, d' \in \C_0, i(e,e) \geq 1, e \mbox{ a geodesic arc}\}
 \]
 be the set of non-simple concatenations of simple arcs. Moreover, let
 \[
   \C_2 = \{ f = e e' \ | \ e,e' \in \C_1, f \mbox{ a geodesic arc}\}
 \]
 Thus, each arc $f = e e' \in \C_2$ has at least two self-intersections, one from $e$ and one from $e'$. Let $\C_1(L)$ and $\C_2(L)$ be the arcs in $\C_1$ and $\C_2$, respectively, that have length at most $L$. For our $\gamma$, let $f_i = e_i e_{i+2}$. The arcs $f_1, \dots, f_m$ are well-defined as long as $m \geq 4$. If $m < 4$, then the lemma holds for any constant $c \geq 3$, because $K \geq 1$.

 It turns out to be easier to bound $m$ using the non-simple arcs $f_1, \dots, f_m$ rather than the simple arcs $\delta_1, \dots, \delta_m$. We will show the following:
 \begin{itemize}
  \item Let $\F = \{ f_1, \dots, f_m\}$. We exhaust $\F$ by sets $\F_1, \dots, \F_N$, where each $\F_i$ is a maximal subset of pairwise disjoint arcs. We bound the size of $\F$ by bounding the size of each $\F_i$, and then by bounding their number, $N$.
  
  \item We show that
  \[
   \# \F_i \leq 2g-2
  \]
  for each $i$ (Lemma \ref{lem:MaxDisjoint}.) We do this by assigning each $f \in \F_i$ to either a pair of pants or a torus with one boundary inside $\S$. Then we show that the set of pairs of pants and one-holed tori assigned to $\F_i$ are all distinct and are part of a pants decomposition of $\S$. This is where we use that each $f \in \C_2$ has at least two self-intersections.
  
  \item Now we want to show that the number $N$ of the sets $\F_1, \dots, \F_N$ satisfies $N \leq c' \sqrt K$ for some universal constant $c'$. (This is proven as part of Lemma \ref{lem:NumberOfDeltaInK}.) We have that 
  \[
   \sum_{i,j = 1}^m i(f_i, f_j) \lesssim i(\gamma, \gamma)
  \]
 where $A \lesssim B$ if $ A \leq c B$ for some universal constant $c$. If each term $i(f_i, f_j)$ contributed at least 1 to this sum, we would be done. Unfortunately, this is not the case, precisely because we can find maximal disjoint subsets $\F_1, \dots, \F_N$ of $\F$. Fortunately, we can also say
 \[
  \sum i(\F_i, \F_j) \lesssim i(\gamma, \gamma)
 \]
 This is good because $i(\F_i, \F_j) \geq 1$ for each $i, j$, by the maximality of these sets. Therefore,
 \[
  \sum_{i, j = 1}^N 1 \lesssim i(\gamma, \gamma)
 \]
 which implies that 
 \[
  N \leq c' \sqrt K
 \]
 for some universal constant $c'$.

 \item Combining the above two statements allows us to show that $\# F \leq c \sqrt K$ for constant $c$ depending only on $\S$ (Lemma \ref{lem:NumberOfDeltaInK}).
 \item Lastly, we give a quick proof that if $\gamma = \delta_1, \dots, \delta_n$, then $n \leq \frac{L}{l_0}$, where $l_0$ is the length of the shortest closed geodesic of $X_0$ to complete the proof (Section \ref{sec:NumDelta}).
 \end{itemize}

\subsection{Maximal sets of pairwise disjoint, non-simple arcs}

The following lemma tells us that a set of pairwise disjoint curves from $\C_2$ cannot have very many elements.

\begin{lem}
\label{lem:MaxDisjoint}
 Fix a set $\{f_1, \dots, f_m\}$ of arcs in $\C_2$.  Suppose $i(f_i, f_j) = 0, \forall i \neq j$. Then $m \leq 2g-2$.
\end{lem}

\begin{proof}
 Because each $f_i$ has at least 2 self-intersections, there is a sense in which it fills either a pair of pants or one-holed torus $\p_i$ (Claims \ref{cla:FindingEights} and \ref{cla:EssentialPants}). We then show that if $\p_i$ and $\p_j$ have an essential overlap as subsurfaces of $\S$, then $f_i$ and $f_j$ must intersect (Claim \ref{cla:EssentialOverlap}). From this we deduce that the set $\p_1, \dots, \p_m$ associated to $f_1, \dots, f_m$ must be part of a pants decomposition of $\S$, and therefore $m \leq 2g-2$.
 
 We start by choosing an $\epsilon$ and arcs $(f_i)_\epsilon \sim_\epsilon f_i$ that we will use in all the claims used to prove this lemma. Fix an $L$ so that $l_0(f_i) \leq L$ for each $i$. Because $\# [\C_0(L) \cup \C_1(L) \cup \C_2(L)] < \infty$, there is some $\eta_L > 0$ so that $\forall \epsilon < \eta_L$, $\forall \delta \in \C_0(L) \cup \C_1(L) \cup \C_2(L)$ and $\forall \delta_\epsilon \sim_\epsilon \delta$, we have that
 \[
  i(\delta,\delta) \leq \# \delta_\epsilon \cap \delta_\epsilon
 \]
 In other words, this choice of $\epsilon$ guarantees that any $\epsilon$-homotopy of any arc $\delta$ we consider will see all of the self-intersections of $\delta$.
 So choose
 \[
  \epsilon < \min\{ \epsilon_L, \eta_L\}
 \]
 where $\epsilon_L$ is the constant from Lemma \ref{lem:GEpsilonDEpsilon}.
 
 For each $i$, choose arcs $(f_i)_\epsilon \sim_\epsilon f_i$ that are geodesics in a negatively curved metric, and so that for each $i$ and $j$, $i(f_i, f_j) = \# (f_i)_\epsilon \cap (f_j)_\epsilon$. Note that $i(f_i, f_j)= 0, \forall i \neq j$ implies that $(f_1)_\epsilon, \dots, (f_m)_\epsilon$ are pairwise disjoint. Suppose each $(f_i)_\epsilon$ is parameterized as an arc, $(f_i)_\epsilon : [0,1] \rightarrow \S$. We find the $\p_i$ using the following, rather technical, claim.
 \begin{claim}
 \label{cla:FindingEights}
  For each $i$, there is a closed sub-interval $I_i \subset [0,1]$ so that
  \begin{itemize}
   \item $(f_i)_\epsilon |_{I_i^0}$ is a simple arc
   \item $(f_i)_\epsilon (\partial I_i) \subset(f_i)_\epsilon (I_i^0)$
  \end{itemize}
 In other words, each $(f_i)_\epsilon$ has a subarc with exactly two self-intersections, which looks like a figure eight. (See Figure \ref{fig:FigureEight}.)
 \end{claim}
 
 \begin{proof}
 For each $f_i$, there are some arcs $e_i, e_i' \in \C_1$ so that $f_i = e_i e_i'$. Thus, we can find subarcs $(e_i)_\epsilon$ and $(e_{i+1})_\epsilon$ of $(f_i)_\epsilon$ with disjoint domains so that $(e_i)_\epsilon \sim_\epsilon e_i$, $(e_i')_\epsilon \sim_\epsilon e_i'$. Because of our choice of $\epsilon$,
 \begin{align*}
 \#(e_i)_\epsilon \cap (e_i)_\epsilon &\geq i(e_i, e_i) \geq 1 \mbox{ and} \\
  \#(e_i')_\epsilon \cap (e_i')_\epsilon &\geq i(e_i', e_i') \geq 1
 \end{align*}
 Because $(e_i)_\epsilon$ has at least one self-intersection point, it has a subarc $\alpha_i$ that can be closed up into a simple closed curve (See Figure \ref{fig:SubLoop}.) That is, there are some $t_i < s_i \in [0,1]$ so that
 \begin{itemize}
  \item $(f_i)_\epsilon |_{(t_i, s_i)}$ is simple
  \item $(f_i)_\epsilon(t_i) = (f_i)_\epsilon (s_i)$
 \end{itemize}

 \begin{figure}[h!] \centering
  \includegraphics{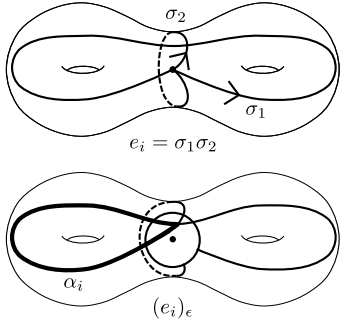}
  \caption{$i(e_i, e_i) > 0$, so we can find a simple subloop $\alpha_i \subset (e_i)_\epsilon$.}
  \label{fig:SubLoop}
 \end{figure}
 
 Let
 \[
  r_i = \min\{ r > s_i \ |\ \exists q, t_i < q < r, (f_i)_\epsilon(r) = (f_i)_\epsilon (q)\}
 \]
In other words, $r_i$ is the first time after $s_i$ that the arc starting at $t_i$ loops back on itself.

Note that $r_i$ exists. If it did not, then $(f_i)_\epsilon$ would be simple on the interval $(t_i, 1]$. But we know that $(e_i')_\epsilon$ is non-simple and its domain comes after the domain of $(e_i)_\epsilon$. 

Let $I_i = [t_i, r_i]$. Then $(f_i)_\epsilon (t_i) = (f_i)_\epsilon (s_i)$ and $(f_i)_\epsilon (r_i)  = (f_i)_\epsilon (q_i)$ for some $s_i, q_i \in (t_i, r_i)$. Furthermore, since we chose $r_i$ to be minimal, $(f_i)_\epsilon$ must be simple on $(t_i, r_i)$. (See Figure \ref{fig:FigureEight}.)

 \begin{figure}[h!] \centering
  \includegraphics{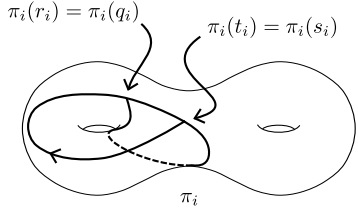}
  \caption[The figure eight subarc of $f_i$]{The figure eight subarc of $f_i$, denoted $\pi_i$.}
  \label{fig:FigureEight}
 \end{figure}
 \end{proof}
 
 For each $i$, let 
 \[
 \pi_i = (f_i)_\epsilon|_{I_i}
 \]
 for the interval $I_i$ from Claim \ref{cla:FindingEights}. Let $N(\pi_i)$ be a regular neighborhood of the graph of $\pi_i$ in $\S$.
 \begin{claim}
 \label{cla:EssentialPants}
   $N(\pi_i)$ is either a pair of pants or one-holed torus, and $\partial N(\pi_i)$ is a set of essential curves.
 \end{claim}
 
 \noindent \textbf{NB:} The curve $\pi_i$ in Figure \ref{fig:FigureEight} fills a one-holed torus, while the curve $\pi_i$ in Figure \ref{fig:NofEight} fills a pair of pants.

 \begin{proof}
 By looking at the Euler characteristic of the graph of $\pi_i$, we see that $N(\pi_i)$ is either a pair of pants or a one-holed torus. To simplify notation, let $N(\pi_i) = \p_i$. We just need to show that no component of $\partial \p_i$ is null-homotopic.
 
 There are two cases. First suppose that $\p_i$ is a one-holed torus with boundary curve $a$. If $a$ is null-homotopic then it bounds a disc $D$. This disc cannot lie in $\p_i$. So if we sew $D$ onto $\p_i$ at $a$, we will get a closed torus. But $\S$ is connected, and it is not a closed torus, so this is a contradiction.
 
 Now suppose that $\p_i$ is a pair of pants. For what follows, refer to Figure \ref{fig:NofEight}.
 
  \begin{figure}[h!] \centering
  \includegraphics{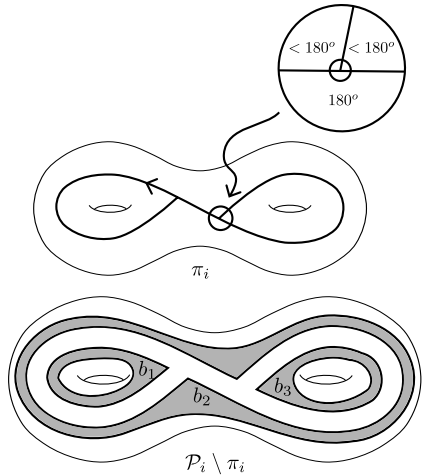}
  \caption{The case when $\p_i$ is a pair of pants}
  \label{fig:NofEight}
 \end{figure}

 Recall that we chose $(f_i)_\epsilon$ to be a geodesic arc in some negatively curved metric $X_\epsilon$. Take any point $x$ on the graph of $\pi_i$, and take a small disk $D$ around it. If we cut $D$ along $\pi_i$, one of two things can happen. If we chose one of the two endpoints of $\pi_i$, then $x$ is a point of intersection between an end of $\pi_i$ and a two-sided subarc of $\pi_i$. In this case, $D \setminus \pi_i$ has three components. The point $x$ lies on the boundaries of these components. In one of the components, the angle at $x$ is $180^o$ with respect to $X_\epsilon$, and in the other two, the angle is strictly smaller than $180^o$. If we choose any other point $x$, then $D \setminus \pi_i$ has two components that have $x$ on their boundaries, and the angle at $x$ is exactly $180^o$ on both of them.
 
 Take $\p_i \setminus \pi_i$. Take the closure of each component of this set separately, and consider the disjoint union of these components. Abusing notation, we will still call the result $\p_i \setminus \pi_i$. Then $\p_i \setminus \pi_i$ has three cylindrical components $C_1, C_2$ and $C_3$. They each have boundary components, denoted $b_1, b_2$ and $b_3$, respectively, that lie on the graph of $\pi_i$. There are exactly four points $x_1, \dots, x_4$ on $b_1 \cup b_2 \cup b_3$ where the angle inside $\p_i \setminus \pi_i$ at $x_j$ is smaller than $180^o$. 
 
Each of $b_1, b_2$ and $b_3$ must contain one of $x_1, \dots, x_4$. Otherwise, they would be simple closed geodesics. But this is impossible since the graph of $\pi_i$ does not contain a simple closed geodesic, as $\pi_i$ is not a simple closed geodesic itself. Thus, without loss of generality, $b_1$ contains $x_1$, $b_2$ contains $x_2$ and $b_3$ contains $x_3$ and $x_4$. Therefore, for each $j$, the sum of exterior angles around $b_j$ is strictly smaller than $360^o$. But if $b_j$ were null-homotopic, it would bound a disc. By Gauss-Bonnet, the sum of exterior angles about the boundary of a disc in the negatively curved metric $X_\epsilon$ is greater than $360^o$. Thus, $b_1$, $b_2$ and $b_3$ are not null-homotopic.

 Therefore, $\p_i$ is either an essential one-holed torus or an essential pair of pants embedded in $\S$.
 
 \end{proof}
 
 From now on, we denote the neighborhood $N(\pi_i)$ by $\p_i$, for each $i$. So we have assigned each $f_i$ an essential pair of pants or one-holed torus $\p_i$. We want to show that $\p_1, \dots, \p_m$ are distinct and that they are part of a pants decomposition of $\S$, which will imply that $m \leq 2g-2$.
 
 \begin{defi}
  Given two sub-surfaces $\S_1, \S_2 \subset \S$, we say $i(\S_1, \S_2) \neq 0$ if for any subsurfaces $\S_1'$ and $\S_2'$ isotopic to $\S_1$ and $\S_2$, respectively, $\S_1' \cap \S_2' \neq \emptyset$. Otherwise, we say $i(\S_1, \S_2) = 0$.
 \end{defi}

 The following claim tells us that if $(f_1)_\epsilon, \dots, (f_m)_\epsilon$ are all pairwise disjoint, then $\p_1, \dots, \p_m$ can all be realized disjointly.
 \begin{claim}
 \label{cla:EssentialOverlap}
  If $i(\p_i, \p_j) \neq 0$, then $\# \pi_i \cap \pi_j \geq 1$.
 \end{claim}
 \begin{proof}
 Suppose $i(\p_i, \p_j) \neq 0$. We will use the following fact to find intersections between $\pi_i$ and $\pi_j$: Because $\p_i = N(\pi_i)$ is a regular neighborhood of $\pi_i$, there is a deformation retract of $\p_i$ onto $\pi_i$. Thus, any closed curve in $\p_i$ is freely homotopic to a closed curve whose image lies in the graph of $\pi_i$. Consider $\partial \p_i$ and $\partial \p_j$ as multicurves. There are two cases to consider. 
 
 The first case is when $i(\partial \p_i, \partial \p_j) \neq 0$ (Figure \ref{fig:PantsBoundaryOverlap}.) In this case, there are some boundary components $a \subset \partial \p_i$ and $b \subset \partial \p_j$ with $i(a,b) \neq 0$. Take closed curves $a'$ and $b'$ whose images in $\S$ lie inside the graphs of $\pi_i$ and $\pi_j$, so that $a'$ is freely homotopic to $a$ and $b'$ is freely homotopic to $b$. By definition, $i(a,b) \neq 0$ implies that $a' \cap b' \neq \emptyset$. But then, $\pi_i \cap \pi_j \neq \emptyset$ as well. 
 
 \begin{figure}[h!] \centering
  \includegraphics{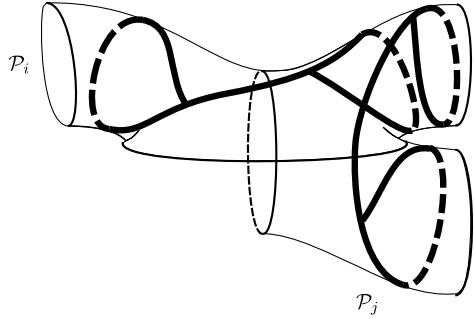}
  \caption{The case when $\partial \p_i \cap \partial \p_j \neq \emptyset$.}
  \label{fig:PantsBoundaryOverlap}
 \end{figure}

 Now consider the case when $i(\partial \p_i, \partial \p_j) = 0$ (Figure \ref{fig:PantsNoBoundaryOverlap}.) Because $\p_i$ and $\p_j$ have an essential overlap as subsurfaces of $\S$, and because they are either pairs of pants or one-holed tori, we can say without loss of generality that  $\p_i$ is isotopic to a subsurface of $\p_j$. (If $\p_i$ and $\p_j$ are both pairs of pants, or if they are both one-holed tori, then they would be isotopic. But if $\p_i$ is a pair of pants and $\p_j$ is a one-holed torus, then it is the closure of $\p_i$ inside $\S$ that is isotopic to $\p_j$.)
 
  \begin{figure}[h!] \centering
  \includegraphics{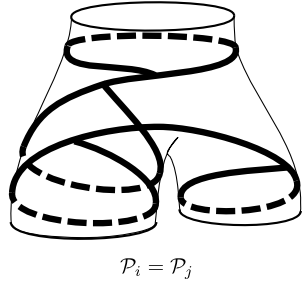}
  \caption{The case when $i(\partial \p_i, \partial \p_j) = 0$.}
  \label{fig:PantsNoBoundaryOverlap}
 \end{figure}
 
 So we can choose a curve $\eta$ with $i (\eta, \eta) = 1$ that can be isotoped to lie inside $\p_i$ and $\p_j$. Let $\eta_i$ and $\eta_j$ be the curves freely homotopic to $\eta$ that lie in the graphs of $\pi_i$ and $\pi_j$, respectively. As $i(\eta,\eta) = 1$, we have that $\# \eta_i \cap \eta_j \geq 1$. So $\# \pi_i \cap \pi_j \geq 1$.
 
 \end{proof}
 
 We have that $i(f_i, f_j) = 0, \forall i \neq j$. We chose $(f_1)_\epsilon, \dots, (f_m)_\epsilon$ so that $\# (f_i)_\epsilon \cap (f_j)_\epsilon = i(f_i, f_j) = 0, \forall i \neq j$. For each $i$, $(f_i)_\epsilon$ corresponds to some pair of pants or one-holed torus $\p_i$ (Claim \ref{cla:EssentialPants}.) 
 We have that $i(\p_i, \p_j) = 0$ for all $i \neq j$ 
 because $(f_i)_\epsilon$ and $(f_j)_\epsilon$ are pairwise disjoint
 (Claim \ref{cla:EssentialOverlap}.)  Thus, $\p_1, \dots, \p_m$ are part of some pants decomposition of $\S$. Therefore, $m \leq 2g-2$.
 
\end{proof}

 \subsection{Bound on number of simple arcs in terms of intersection number}
We are finally ready to bound the number of simple arcs in a geodesic $\gamma \in \G_*^c$ in terms of its self-intersection number.
\begin{lem}
\label{lem:NumberOfDeltaInK}
 Let $\gamma \in \G_*^c(L,K)$. Suppose
 \[
  \gamma = \delta_1 \dots \delta_m
 \]
is the shortest way to write $\gamma$ as a concatenation of arcs $\delta_1, \dots, \delta_m \in \C_0$. If $K \geq 1$, then
\[
 m \leq c \sqrt K
\]
where $c$ is a constant depending only on the topology of $\S$.
\end{lem}

\begin{proof}
Suppose $\gamma = \delta_1 \dots \delta_m$, where $\delta_i \in \C_0, \forall i$. If $m \leq 3$, then we are done for $c \geq 3$. So suppose $m \geq 4$. Let $e_i = \delta_i \delta_{i+1}$. Then $e_i \in \C_1$. Now let $f_i = e_i e_{i+2}$. Then $f_i \in \C_2$. Both $e_i$ and $f_i$ are well-defined for each $i$ because there are at least four distinct simple arcs in $\gamma$.

First we show that
\[
 \sum_{i, j = 1}^m i(f_i, f_j) \leq 32 i(\gamma, \gamma)
\]
Take the $\epsilon$ we defined at the start of the proof of Lemma \ref{lem:MaxDisjoint}. Take a curve $\gamma_\epsilon \sim_\epsilon \gamma$ so that $i(\gamma, \gamma) = \# \gamma_\epsilon \cap \gamma_\epsilon$. 

The homotopy from $\gamma$ to $\gamma_\epsilon$ gives a correspondence between $f_i$ and some curve $(f_i)_\epsilon \subset \gamma_\epsilon$ for each $i$. If we parameterize $\gamma: [0,1] \rightarrow \S$ and let $J_1, \dots, J_m$ be the domains of $f_1, \dots, f_m$, respectively, then each $t \in [0,1]$ lies in exactly 4 of $J_1, \dots, J_m$. Thus, exactly 16 pairs $(f_i)_\epsilon, (f_j)_\epsilon$ intersect at each self-intersection point of $\gamma_\epsilon$. Because $\epsilon < \eta_L$, $i(f_i, f_j) \leq \#(f_i)_\epsilon \cap (f_j)_\epsilon$. Thus,
\[
 \frac 12 \sum_{i,j} i(f_i, f_j) \leq 16 i(\gamma, \gamma)
\]
where the $\frac 12$ comes from the fact that each pair $f_i, f_j$ appears twice in the sum on the left hand side.

Let $\F = \{f_1, \dots, f_m\}$. We want to show $i(f_i, f_j) \neq 0$ for sufficiently many pairs $i,j$. To do this, we decompose $\F$ into small sets $\F_1, \dots, \F_N$ for which $i(\F_i, \F_j) > 0, \forall i,j$. Up to renumbering indices, let 
\[
 \F_1 = \{f_1, \dots, f_{n_1}\}
\]
be a maximal subset of $\F$ so that $i(f_i, f_j) = 0, \forall i \neq j = 1, \dots, n_1$. Given $\F_1, \dots, \F_i$, let 
\[
 \F_{i+1} = \{f_{n_i+1}, \dots, f_{n_{i+1}}\}
\]
be a maximal subset of $\F \setminus (\F_1 \cup \dots \cup \F_i)$ so that $i(f_i, f_j) = 0, \forall i \neq j = n_i+1, \dots,n_{i+1}$. Again, this is up to renumbering. In this way, we exhaust $F$ with a list $\F_1, \dots, \F_N$ of such subsets.

Whenever $f_i, f_j$ are in the same set $\F_k$, $i(f_i, f_j) = 0$. Thus,
\[
 \sum_{ i \geq j} i(\F_i, \F_j) = \sum_{i \geq j} i(f_i, f_j)
\]
Thus, 
\[
 \sum_{i \geq j} i(\F_i, \F_j) \leq 16 i(\gamma, \gamma)
\]
Because each $\F_i$ is maximal in what is left over when we take away all previous sets, each $f \in \F_j$ intersects some $f' \in \F_j$, when $i > j$. In other words, 
\[
 i(\F_i, \F_j) \neq 0
\]
where this intersection number is the sum of intersections of each element in $\F_i$ and each element in $\F_j$.

Combining all of this, we get that 
\[
 \frac{N(N-1)}{2} \leq 16 i(\gamma, \gamma)
\]
%
%
%
%
%
%
%
%
Thus, for some universal constant $c'$, we get that
\[
 N \leq c' \sqrt K
\]
(For example $c'$ can be taken smaller than 15.)

We wish to bound $m = \# \F$. But
\[
 \#\F = \# \F_1 + \dots + \# \F_N
\]
 The sets $\F_1, \dots, \F_N$ satisfy the conditions of Lemma \ref{lem:MaxDisjoint}, so
\[
 \#\F_i \leq 2g-2
\]
for each $i$. 
So we get
\[
 m \leq c \sqrt K
\]
for $c = (2g-2)c'$. Note that $c$ depends only on the topology of $\S$.
\end{proof}
 
\subsection{Bound on number of simple arcs in terms of length and intersection number}
\label{sec:NumDelta}
We are done with the proof of Lemma \ref{lem:DeltaLengthCount} as soon as we bound the length of the sequence $\gamma = \delta_1 \dots \delta_m$ in terms of $l_0(\gamma)$. Each $\delta_i$ consists of at least one saddle connection. Let $l_0$ be the length of the shortest closed geodesic on $X_0$. Then $l_0(\gamma) \geq l_0 m$. If $l_0(\gamma) \leq L$, then
\[
 m \leq \frac{L}{l_0}
\]

So, we have shown that
\[
 m \leq \min \{\frac{L}{l_0},c \sqrt K\}
\]
where $c$ depends just on $\S$, and $l_0$ depends on $X_0$.

\end{proof}

\section{Size of $\G_*^c(L,K)$}
We have shown that any $\gamma \in \G_*^c(L,K)$ can be written as $\gamma = \delta_1 \dots \delta_m$ for $\delta_i \in \C_0$ for each $i$, and $m \leq \min \{\frac{L}{l_0},c \sqrt K\}$ (Lemma \ref{lem:DeltaLengthCount}). Furthermore, if $l_0(\gamma) \leq L$, then $l_0(\delta_i) \leq L$ for each $i$. We have shown that $\# \C_0(L) \leq c_0 L^{c_g}$ where $c_0$ depends only on $X_0$ and $c_g$ depends only on $\S$ (Lemma \ref{lem:SimpleCount}). The number of possible sequences $\delta_1, \dots, \delta_n$ that satisfy these criteria is at most
\[
 (c_0 L^{c_g})^{\min \{\frac{L}{l_0},c \sqrt K\}}
\]

The map from geodesics $\gamma$ to ordered sets $\delta_1, \dots, \delta_m$ is injective as any sequence $\delta_1, \dots, \delta_m$ determines a closed curve and geodesics in $\G_*$ are unique in their free homotopy classes. So the number of all possible sequences $\delta_1, \dots, \delta_m$ for which $m \leq \min \{\frac{L}{l_0},c \sqrt K\}$ and $l_0(\delta_i) \leq L, \forall i$ bounds $\# \G_*^c(L,K)$ from above. So in particular,
\[
 \# \G_*^c(L.K) \leq
  (c_0 L^{c_g})^{c \sqrt K}
\]

 \bibliographystyle{alpha}
  \bibliography{recount,researchForRecount}
\end{document}